\documentclass[%
 aip,
sd,%
 amsmath,amssymb,
reprint,%
]{revtex4-1}
\pdfoutput=1
\usepackage{graphicx}
\usepackage{dcolumn}
\usepackage{bm}
\usepackage{graphicx}
\usepackage{epstopdf}
\usepackage{subcaption}
\usepackage{color}
\usepackage{rotating}
\newtheorem{theorem}{Theorem}
\newcommand{\change}[1]{{\color{black}#1}}

\draft 

\begin{document}


\title{Early-warning indicators in the dynamic regime} 



\author{Paul Ritchie}
\email[]{pdlr201@exeter.ac.uk}

\author{Jan Sieber}
\email[]{J.Sieber@exeter.ac.uk}
\affiliation{Centre for Systems, Dynamics and Control, College of Engineering, Mathematics and Physical Sciences, Harrison Building, University of Exeter, Exeter, EX4 4QF, United Kingdom}


\date{\today}

\begin{abstract}
Early-warning indicators (increase of autocorrelation and variance) are commonly applied to time series data to try and detect tipping points of real-world systems. The theory behind these indicators originates from approximating the fluctuations around an equilibrium observed in time series data by a linear stationary (Ornstein-Uhlenbeck) process. Then for the approach of a bifurcation-type tipping point the formulas for the autocorrelation and variance of an Ornstein-Uhlenbeck process detect the phenomenon `critical slowing down'. The assumption of stationarity and linearity introduces two sources of error in the early-warning indicators. We investigate the difference between the theoretical and observed values for the early-warning indicators for the saddle-node normal form bifurcation with linear drift.   
\end{abstract}

\keywords{}

\maketitle 


This paper will focus on two common early-warning indicators, the increase of autocorrelation and variance in time series generated by a dynamical system in which the system parameter slowly approaches a saddle-node bifurcation (disturbed by white noise). Section \ref{sec: Numerical FP} outlines the numerical approach for solving the Fokker-Planck equation and how this is used to calculate the early-warning indicators.

Fluctuations around a quasi-static equilibrium observed in a time series can be approximated linearly by the Ornstein-Uhlenbeck process. For an Ornstein-Uhlenbeck process the autocorrelation and variance is known. The early-warning indicators for the Ornstein-Uhlenbeck process assume quasi-stationarity and linearity and will be referred to as the linear quasi-static indicators. We analyse the systematic differences between the nonlinear dynamic and linear quasi-static early-warning indicators for the saddle-node normal form with linear drift in Section \ref{subsec: lin drift}. Section \ref{subsec: SN nlin drift} compares these early-warning indicators for a nonlinear drift motivated from a model for rate-induced tipping.

We end the paper in Section \ref{sec: ISM} with a case study from climate science, a conceptual model for the Indian summer monsoon. The model developed by \citet{zickfeld2004modeling} shows that an increase of the planetary albedo or decrease in $CO_{2}$ can lead to the summer monsoon being stopped. Simulating a slow increase of planetary albedo such that the system passes through the saddle-node bifurcation we show that the early-warning indicators are present for this tipping event.

\section{Numerical solution of the Fokker-Planck equation}
\label{sec: Numerical FP}

In this section we will discuss the method for calculating the probability density $P(x,t)$ for a random variable $X_{t}$ governed by the SDE:

\begin{equation}
\mathrm{d}X_{t} = f(X_{t},t)\mathrm{d}t + \sqrt{2D}\mathrm{d}W_{t}
\label{sec:SNexamples SDE}
\end{equation}

\noindent where the drift $f(x,t)$ is related to a potential $U(x,t)$ by $f(x,t) = -\partial_{x}U(x,t)$. The time evolution for the probability density $P(x,t)$ of the random variable $X_{t}$ is described by the Fokker-Planck equation:

\begin{equation}
\dfrac{\partial P}{\partial t} = D\dfrac{\partial^{2}P(x,t)}{\partial x^{2}} - \dfrac{\partial}{\partial x}\bigg(f(x,t)P(x,t)\bigg)
\label{Fokker Planck}
\end{equation}

\noindent where we start at $t = t_{0}$ with an initial condition of \eqref{sec:SNexamples SDE} distributed according to some density $P(x,t_{0})$. We proceed by describing how we solve the Fokker-Planck equation before determining the methods for calculating the early-warning indicators and escape rate.

\subsection{Discretisation of equilibrium problem}

Initially we will calculate the stationary solution, $P_{s}(x)$ of the Fokker-Planck equation

\begin{equation}
0 = D\dfrac{\partial^{2}P_{s}(x)}{\partial x^{2}} - \dfrac{\partial}{\partial x}\bigg(f(x)P_{s}(x)\bigg)
\label{stationary FP}
\end{equation}

\noindent on the interval $[x_{\mathrm{start}},x_{\mathrm{end}}]$ with Dirichlet boundary conditions

\begin{equation*}
P_{s}(x_{\mathrm{start}}) = P_{s}(x_{\mathrm{end}}) = 0
\end{equation*}

\noindent Choosing Dirichlet boundary conditions enables us to monitor the escape rate per unit time when we solve the time dependent Fokker-Planck equation \eqref{Fokker Planck}. However, this will mean that we can only obtain an approximation for the stationary solution on a bounded domain, as will be discussed shortly.

We split the closed $x$ domain into $N$ equal intervals of length $\Delta x = (x_{\mathrm{end}}-x_{\mathrm{start}})/(N+1)$, we denote the interval points as \citep{morton2005numerical}

\begin{equation*}
x_{i} = x_{\mathrm{start}} + (i-1)\Delta x \qquad \mbox{for} \quad i = 1,2,...,N+1
\end{equation*}

\noindent Likewise, for an arbitrary function $g(x)$, for notation simplicity we will denote $g(x_{i}) = g_{i}$. We approximate the derivatives of the function $g(x)$ using finite difference methods \citep{smith1985numerical}:

\begin{equation}
\bigg(\dfrac{\mathrm{d}g}{\mathrm{d}x}\bigg)_{i} = \dfrac{g_{i+1} - g_{i-1}}{2\Delta x} + \mathcal{O}(\Delta x)^{2}
\label{1st deriv x}
\end{equation}  

\noindent and

\begin{equation}
\bigg(\dfrac{\mathrm{d}^{2}g}{\mathrm{d}x^{2}}\bigg)_{i} = \dfrac{g_{i+1} - 2g_{i} + g_{i-1}}{(\Delta x)^{2}} + \mathcal{O}(\Delta x)^{2}
\label{2nd deriv x}
\end{equation} 

\noindent by rearranging the stationary Fokker-Planck equation \eqref{stationary FP} to

\begin{equation}
0 = \bigg(D\dfrac{\partial^{2}}{\partial x^{2}} - \dfrac{\partial f(x)}{\partial x} - f(x)\dfrac{\partial}{\partial x}\bigg)P_{s}(x)
\label{stationary FP matrix}
\end{equation}

\noindent and using the finite difference approximations \eqref{1st deriv x} and \eqref{2nd deriv x} we can express \eqref{stationary FP matrix} as

\begin{equation}
0 = AP_{s}
\label{Stationary matrix}
\end{equation}

\noindent where $A$ is an $(N+1)\times(N+1)$ tridiagonal matrix with the non-zero entries given by

\begin{align*}
A(i,i-1) &= \dfrac{f_{i}}{2\Delta x} + \dfrac{D}{(\Delta x)^{2}} \\       
A(i,i) &= \dfrac{f_{i-1} - f_{i+1}}{2\Delta x} - \dfrac{2D}{(\Delta x)^{2}} \\ 
A(i,i+1) &= -\dfrac{f_{i}}{2\Delta x} + \dfrac{D}{(\Delta x)^{2}}   
\end{align*} 

\noindent for $i = 2,...,N$ and then setting Dirichlet boundary conditions on the probability density function we also have

\begin{align*}
A(1,1) &= 1 \\
A(N+1,N+1) &= 1
\end{align*}

\noindent The stationary probability density $P_{s}(x)$ can be determined by \change{replacing \eqref{Stationary matrix}} with the eigenvalue equation: 

\begin{equation*}
A\textbf{v}_{j} = \gamma_{j}\textbf{v}_{j}
\end{equation*}

\noindent where the $\gamma_{j}$'s are the eigenvalues and $\textbf{v}_{j}$'s the corresponding eigenvectors of the matrix $A$. On an infinite domain with a potential $U(x)\to+\infty$ as $x\to\pm\infty$ there exists an eigenvalue $\gamma_{1} = 0$ \citep{risken2012fokker}. However, as previously noted, on a bounded domain with Dirichlet boundary conditions we can only obtain an approximation for the stationary probability density $P_{s}(x)$. This can be seen by integrating the eigenvalue problem

\begin{equation*}
\gamma_{1}\textbf{v}_{1}(x) = \dfrac{\partial}{\partial x}\bigg(D\dfrac{\partial \textbf{v}_{1}(x)}{\partial x} - f(x)\textbf{v}_{1}(x)\bigg)
\end{equation*}

\noindent over the domain $x\in [x_{\mathrm{start}},x_{\mathrm{end}}]$ and because of the Dirichlet conditions imposed we get

\begin{equation*}
\gamma_{1}\int\limits_{x_{\mathrm{start}}}^{x_{\mathrm{end}}}\textbf{v}_{1}(x)\mathrm{d}x = D\big(\textbf{v}_{1}'(x_{\mathrm{end}})-\textbf{v}_{1}'(x_{\mathrm{start}})\big)
\end{equation*}

\noindent Using that at the boundaries $\textbf{v}_{1}'(x_{\mathrm{start}})\geq 0 \geq \textbf{v}_{1}'(x_{\mathrm{end}})$ demonstrates that the leading eigenvalue $\gamma_{1} \leq 0$. We shall now assume $\gamma_{1} = 0$ and prove by contradiction that $\gamma_{1} < 0$. The general solution of \eqref{stationary FP} is given by

\begin{equation*}
P_{s}(x) = C\int\limits_{x_{0}}^{x}\exp\bigg(\dfrac{U(x')-U(x)}{D}\bigg)\mathrm{d}x' 
\end{equation*}

\noindent The boundary condition $P_{s}(x_{\mathrm{start}}) = 0$ determines $x_{0} = x_{\mathrm{start}}$. Though the other boundary condition $P_{s}(x_{\mathrm{end}}) = 0$ is only satisfied for $C = 0$ which corresponds to $P_{s}(x) = 0$ for all $x$. Therefore we must have $\gamma_{1}<0$ and thus, the corresponding normalised eigenvector $\textbf{v}_{1}$ gives an approximation to the stationary probability density $P_{s}(x)$.  

\subsection{Discretisation of evolution equation}  
\label{sec:Doee}

Solving the time dependent Fokker-Planck equation \eqref{Fokker Planck} we also discretise the time domain $[t_{0},T_{\mathrm{end}}]$ into $M$ equal intervals of length $\Delta t$, and denote:

\begin{equation*}
t_{n} = t_{0} + (n-1)\Delta t \qquad \mbox{for} \quad n = 1,2,...,M+1
\end{equation*}

\noindent and again for an arbitrary function $g(x,t)$ for notational simplicity we will write $g(x_{i},t_{n}) = g_{i}^{n}$.

A numerical technique often used for solving a PDE of the form:

\begin{equation*}
\dfrac{\partial P}{\partial t} = F\bigg(P,f,x,t,\dfrac{\partial P}{\partial x},\dfrac{\partial^{2}P}{\partial x^{2}},\dfrac{\partial f}{\partial x}\bigg)
\end{equation*}

\noindent is the implicit Crank-Nicolson method, given as \citep{quarteroni2014scientific}

\begin{equation}
P_{i}^{n+1} = P_{i}^{n} + \dfrac{1}{2}\Delta t\bigg(F_{i}^{n} + F_{i}^{n+1}\bigg)
\label{Crank Nicolson}
\end{equation}

\noindent Applying the Crank-Nicolson method \eqref{Crank Nicolson} to the Fokker-Planck equation \eqref{Fokker Planck} we can write this as

\begin{equation*}
0 = A_{1}^{n}P^{n} + A_{2}^{n+1}P^{n+1}
\end{equation*}

\noindent where $A_{1}^{n}$ and $A_{2}^{n+1}$ are $(N+1)\times(N+1)$ tridiagonal matrices with the non-zero entries given as follows.

\begin{align}
\nonumber
A_{j}^{k}(i,i-1) &= \dfrac{1}{2}\Delta t\bigg(\dfrac{f_{i}^{k}}{2\Delta x} + \dfrac{D}{(\Delta x)^{2}}\bigg) \\  
\label{time depend mat}   
A_{j}^{k}(i,i) &= (-1)^{j+1} + \dfrac{1}{2}\Delta t\bigg(\dfrac{f_{i-1}^{k} - f_{i+1}^{k}}{2\Delta x} - \dfrac{2D}{(\Delta x)^{2}}\bigg) \\ 
\nonumber
A_{j}^{k}(i,i+1) &= \dfrac{1}{2}\Delta t\bigg(-\dfrac{f_{i}^{k}}{2\Delta x} + \dfrac{D}{(\Delta x)^{2}}\bigg)   
\end{align} 

\noindent for $i = 2,...,N$. Note that the matrices $A_{1}^{n}$ and $A_{2}^{n+1}$ are now time dependent for a time dependent drift, $f$. We still apply Dirichlet boundary conditions to the probability density $P^{n+1}$ and so we have

\begin{align*}
A_{2}^{n+1}(1,1) &= 1 \\
A_{2}^{n+1}(N+1,N+1) &= 1
\end{align*}

\noindent Therefore we are now in a position to calculate the time evolution of the probability densities starting with the stationary distribution $P_{s}(x,t_{0})$ using the formula:

\begin{equation}
P^{n+1} = -(A_{2}^{n+1})^{-1}A_{1}^{n}P^{n}
\label{FP sol step}
\end{equation}

\noindent The Fokker-Planck equation \eqref{Fokker Planck} is of the form of an advection-diffusion equation. Though when solving numerically an advection-diffusion equation using the Crank-Nicolson method \eqref{Crank Nicolson} in conjunction with the finite difference methods \eqref{1st deriv x}--\eqref{2nd deriv x} there are restrictions on the step sizes $\Delta x$ and $\Delta t$ which need to be adhered to. 

\subsection{Accuracy and stability conditions}

The finite difference methods used, in particular \eqref{2nd deriv x}, are susceptible to producing spatial oscillations in the tails of the probability density. Though increasing the diffusion coefficient $D$ reduces the oscillations until they finally disappear. This gives rise to a condition on the Peclet number $P_{e}$, which ensures the drift terms are not too large compared to the diffusion terms inside the matrices $A_{1}$ and $A_{2}$ \citep{szymkiewicz2010numerical}:

\begin{equation}
P_{e} = \dfrac{\max|f|\Delta x}{D} \leq 2
\label{Peclet number}
\end{equation}

\noindent Although the Crank-Nicolson scheme is unconditionally stable we need a tight restriction on $\Delta t$ relative to $\Delta x$ to prevent unrealistic behaviour (time dependent oscillations). This is given by \citep{versteeg2007introduction}:

\begin{equation}
\Delta t < \dfrac{(\Delta x)^{2}}{D}
\label{CN stability}
\end{equation} 

\noindent The limitations of this method can be seen from \eqref{Peclet number} and \eqref{CN stability} if we want to consider examples with either a large drift term $f(x,t)$ or a small noise level $D$. We would need to choose the grid spacing $\Delta x$ small to satisfy \eqref{Peclet number}, this affects the time step $\Delta t$, which will also have to be small to satisfy \eqref{CN stability}, and thus increases computational effort.

\subsection{Example: Fokker-Planck equation for straight drift}

\change{We calculate the probability density numerically \eqref{FP sol step} from the Fokker-Planck equation \eqref{Fokker Planck} for an example where the explicit solution is also known.} In this example, we consider the SDE as:

\begin{equation}
\mathrm{d}X_{t} = -\mathrm{d}t + \sqrt{2D}\mathrm{d}W_{t}
\label{FP example SDE}
\end{equation} 

\noindent and thus, $f(x,t) = f = -1$. \change{If we make the change of coordinates $X_{t} = Y_{t} - t$, \eqref{FP example SDE} is transformed into an SDE for pure diffusion:}

\begin{equation}
\mathrm{d}Y_{t} = \sqrt{2D}\mathrm{d}W_{t}
\label{FP example y SDE}
\end{equation}

\noindent which has an explicit solution for the probability density for an initial Dirac delta distribution at $(y,t) = (y_{0},t_{0})$:

\begin{equation*}
P(y,t) = \dfrac{1}{\sqrt{4\pi D(t-t_{0})}}\exp\bigg(-\dfrac{(y-y_{0})^{2}}{4D(t-t_{0})}\bigg)
\end{equation*}

\noindent \change{We apply the transformation back to obtain the density for a straight drift:}

\begin{equation}
P(x,t) = \dfrac{1}{\sqrt{4\pi D(t-t_{0})}}\exp\bigg(-\dfrac{(x-x_{0}+t)^{2}}{4D(t-t_{0})}\bigg)
\label{density straight drift}
\end{equation} 

\begin{figure}[h!]
        \centering
        \subcaptionbox{\label{FP illustration}}[0.45\linewidth]
                {\includegraphics[scale = 0.3]{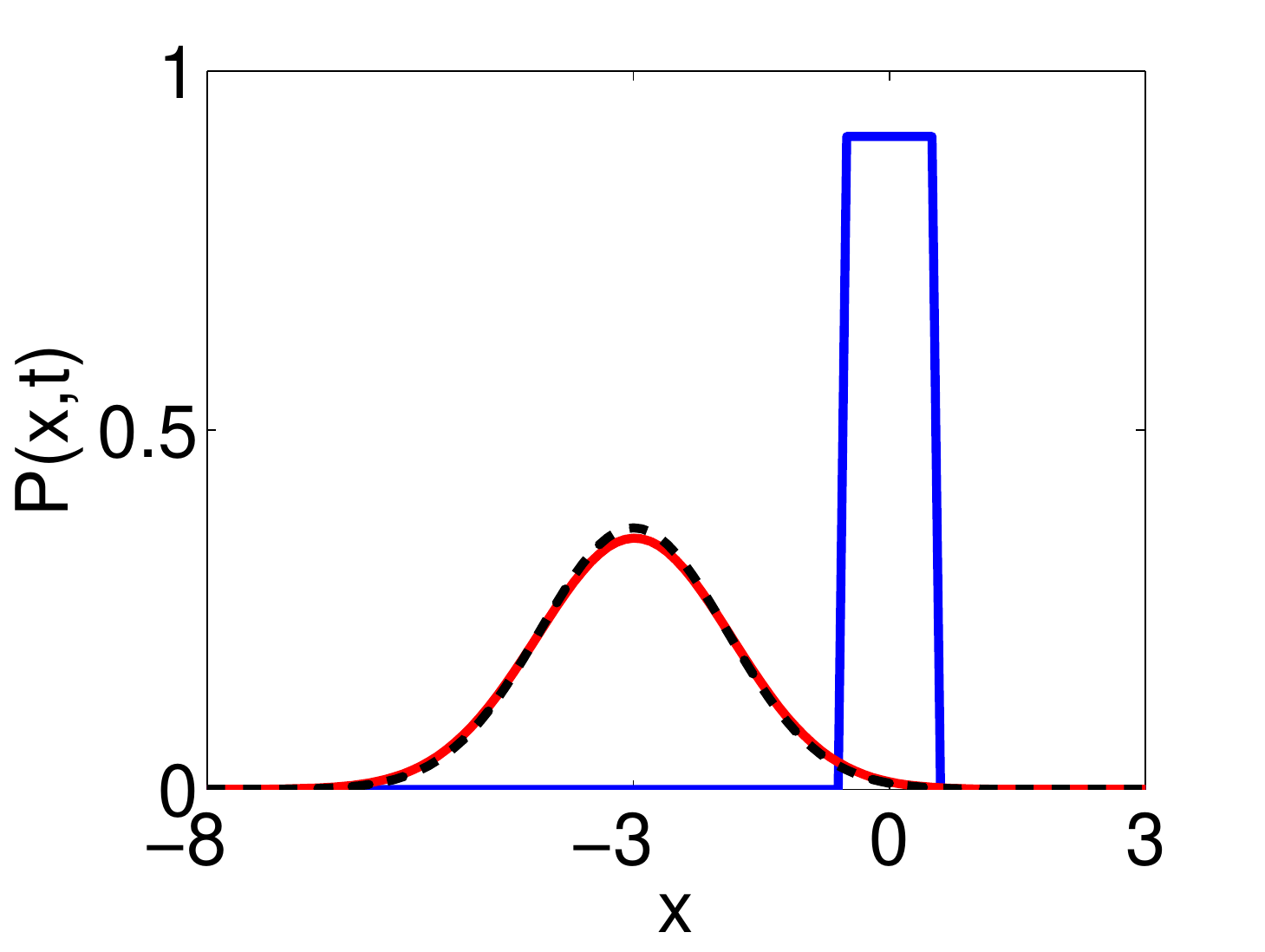}}
       \hfill
        \subcaptionbox{\label{FP illustration density}}[0.45\linewidth]
                {\includegraphics[scale = 0.3]{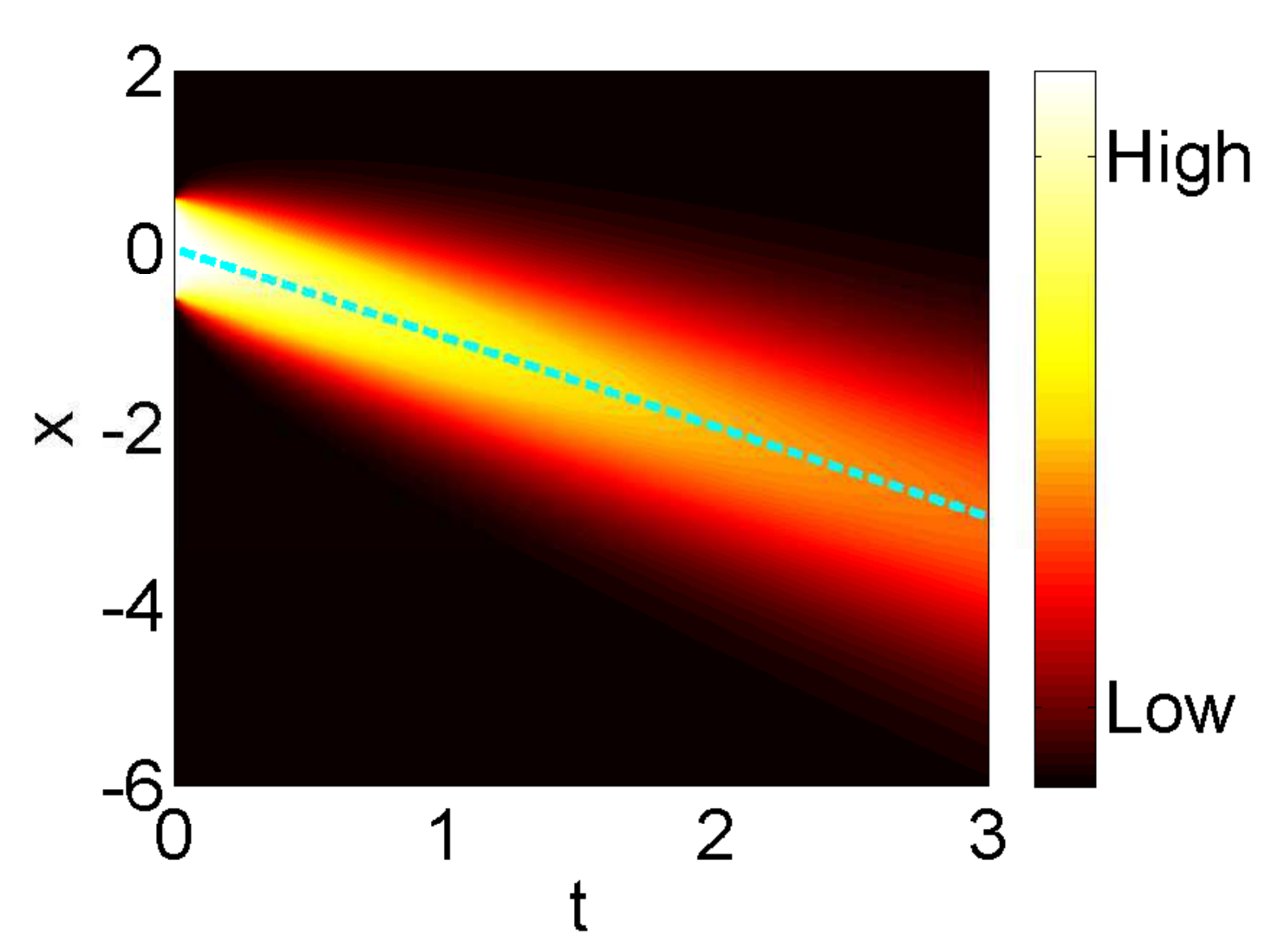}}
        ~ 
        \caption[Illustration of Fokker-Planck equation for straight drift.]{Illustration of the Fokker-Planck equation for straight drift $f = -1$. (a) Numerical initial density $P(x,0)$ (blue) and density at $T_{\mathrm{end}} = 3$ (red). Analytical solution \eqref{density straight drift} for initially starting with Dirac delta distribution at $(x,t_{0}) = (0,0)$ (black dashed). (b) Time space plot of density represented by colour, blue dashed line indicates solution of ODE $\dot{x} = -1$ with initial condition $x_{0} = 0$ at $t_{0} = 0$. Other parameters: $D = 0.2$, $\Delta x = 0.1$, $\Delta t = 0.01$.}\label{FP example}
\end{figure}

\noindent \change{Figure \ref{FP illustration} presents the numerically calculated probability density (red) and the analytical solution, \eqref{density straight drift} (black dashed) after $3$ time units ($t_{0} = 0$, $T_{\mathrm{end}} = 3$). The analytical solution assumes a Dirac delta distribution at $(x_{0},t_{0}) = (0,0)$, whereas the numerical calculation starts from the initial density (blue) given in Figure \ref{FP illustration}. We choose $\Delta x = 0.1$ and $\Delta t = 0.01$ to ensure that both the accuracy \eqref{Peclet number} and stability \eqref{CN stability} conditions are satisfied and thus, reduces oscillations in the tails of the numerical density. We observe there is a good agreement between the analytical and numerically calculated densities. From \eqref{density straight drift} we see that the densities are centred about $x = -T_{\mathrm{end}} = -3$ and the variance of the densities at $T_{\mathrm{end}} = 3$ are proportional to the noise level $D$.}  

Furthermore, Figure \ref{FP illustration density} gives an illustration as to how the numerically calculated probability density spreads out over time indicated by the colour. Initially the density is narrow and tall centred about the deterministic solution ($D = 0$) of \eqref{FP example SDE} with initial condition $(x_{0},t_{0}) = (0,0)$ (blue dashed line). As the time increases the density widens but remains centred about the deterministic solution. 

\subsection{Numerical calculation of early-warning indicators and escape rate}  
\label{subsec: Num EWI}

Early-warning indicators are usually applied to either real-world time series data or time series calculated from an ensemble of realisations using the SDE \eqref{sec:SNexamples SDE}. The indicators are used to quantify changes in the statistical properties of the time series data on the approach to a tipping point. In practice a `sliding window' is used to smooth out the noise in time series data \citep{dakos2012methods}. Though the results of this depend on the amount of time series data available, as this will affect the width of the sliding window. Whereas, the Fokker-Planck equation \eqref{Fokker Planck} gives us the opportunity to compute them directly. Two early-warning indicators that we will calculate numerically for the examples we consider are: increased variance and increased lag-1 autocorrelation. We will now outline the numerical algorithms used to calculate the variance and lag-1 autocorrelation and where the relevant theory can be found in \citet{ross2004introduction}.

We will again work on the closed domain $x \in [x_{\mathrm{start}},x_{\mathrm{end}}]$ discretised into $N$ equal intervals of length $\Delta x$, the interval points given as

\begin{equation*}
x_{i} = x_{\mathrm{start}} + (i-1)\Delta x \qquad \mbox{for} \quad i = 1,2,...,N+1
\end{equation*}

\noindent Likewise, the time domain $[t_{0},T_{\mathrm{end}}]$ is discretised into $M$ equal intervals of length $\Delta t$, where

\begin{equation*}
t_{n} = t_{0} + (n-1)\Delta t \qquad \mbox{for} \quad n = 1,2,...,M+1
\end{equation*}

\noindent For a random variable $X_{t_{n}}$ that is the solution of the SDE \eqref{sec:SNexamples SDE} the probability density function $P(x,t_{n})$ at time $t_{n}$ can be obtained via the methods presented in Section \ref{sec:Doee}.

\paragraph{Variance:} The variance of the stochastic process $X_{t_{n}}$ at time $t_{n}$ is defined as

\begin{equation}
\mathrm{Var}(X_{t_{n}}) = \mathbb{E}(X_{t_{n}}^{2}) - \mathbb{E}(X_{t_{n}})^{2}
\label{sec:SN variance}
\end{equation}

\noindent where for a discrete random variable $X_{t_{n}}$ with probability density $P(x,t_{n})$, the $k^{\mathrm{th}}$ moment is given by:

\begin{equation*}
\mathbb{E}(X_{t_{n}}^{k}) = \sum_{i=1}^{N+1}x_{i}^{k}P(x_{i},t_{n})
\end{equation*}

\noindent for $k \geq 1$.

\paragraph{Lag-1 autocorrelation:} The autocorrelation $a_{n}$ given at time $t_{n}$ is the correlation between the random variables $X_{t_{n-1}}$ and $X_{t_{n}}$, separated by a time step $\Delta t$:

\begin{equation}
a_{n} = \dfrac{\mathrm{Cov}(X_{t_{n-1}},X_{t_{n}})}{\sqrt{\mathrm{Var}(X_{t_{n-1}})\mathrm{Var}(X_{t_{n}})}}
\label{autocorrelation SN sec}
\end{equation}

\noindent We have already shown how to calculate the variance of a random variable and so we will now focus on the covariance between two random variables. To simplify notation, we make the substitutions $X = X_{t_{n-1}}$ and $Y = X_{t_{n}}$, then the covariance between $X$ and $Y$ is given as follows:

\begin{equation*}
\mathrm{Cov}(X,Y) = \mathbb{E}(XY) - \mathbb{E}(X)\mathbb{E}(Y)
\end{equation*}

\noindent The expectation $\mathbb{E}(XY)$ can be calculated provided the joint probability density function of $X$ and $Y$ is known, denoted $P_{X,Y}(x_{i},y_{j})$:

\begin{equation*}
\mathbb{E}(XY) = \sum_{j=1}^{N+1}\sum_{i=1}^{N+1}x_{i}y_{j}P_{X,Y}(x_{i},y_{j})
\end{equation*}

\noindent The joint probability density can be expressed as follows:

\begin{equation}
P_{X,Y}(x_{i},y) = P_{Y|X}(y|x_{i})P(x_{i})
\label{joint probability density}
\end{equation}

\noindent where $P_{Y|X}(y|x_{i})$ is the conditional probability density for $Y$ given $X = x_{i}$ and $P(x_{i})$ is the probability for the random variable $X = x_{i}$.

The conditional probability $P_{Y|X}(y|x_{i})$ in \eqref{joint probability density} is obtained first by setting an approximate Dirac delta distribution $P_{x_{i}}^{0}$ at the point $x_{i}$. For example

\begin{equation*}
P_{x_{i}}^{0}(x_{k}) = \begin{cases} \frac{1}{\Delta x} \qquad &\mbox{if} \quad k=i \\ 0 \qquad &\mbox{if} \quad k\neq i \end{cases}
\end{equation*}

\noindent such that the area under $P_{x_{i}}^{0}$ is approximately equal to $1$. The probability density $P_{x_{i}}^{0}$ represents the "given $X = x_{i}$" part (i.e. $P(X=x_{i}) = 1$) in the conditional probability. Therefore the conditional probability $P_{Y|X}(y|x_{i})$ is calculated by evolving the density $P_{x_{i}}^{0}$ over a single time step $\Delta t$ using \eqref{FP sol step}:

\begin{equation*}
P_{Y|X}(y|x_{i}) = P_{x_{i}}^{1} = -(A_{2}^{n+1})^{-1}A_{1}^{n}P_{x_{i}}^{0} 
\end{equation*}  

\noindent where $A_{j}^{k}$ were defined in equation \eqref{time depend mat}. 

\paragraph{Escape rate:} The escape rate, $r$, evaluates the amount of escape from the domain per unit time and is calculated using the probability density $P(x,t_{n})$. For a probability density $P_{\infty}(x,t_{n})$ on the infinite domain we would have:

\begin{equation*}
\int\limits_{-\infty}^{\infty}P_{\infty}(x,t_{n})\mathrm{d}x = 1 \qquad \forall n
\end{equation*}

\noindent However, our probability density $P(x,t_{n})$ is on a bounded domain, $x \in [x_{\mathrm{start}},x_{\mathrm{end}}]$, with Dirichlet boundary conditions. The probability density $P^{n}(x)$ calculated using \eqref{FP sol step} (with normalised density $P^{n-1}(x) := \tilde{P}^{n-1}(x)$, such that $\int \tilde{P}^{n-1}(x)\mathrm{d}x = 1$) will have a 'survival rate', $s^{n}$, approximated by the trapezoidal rule \citep{larson2009calculus}:

\begin{align*}
s^{n} &= \int\limits_{x_{\mathrm{start}}}^{x_{\mathrm{end}}}P^{n}(x)\mathrm{d}x \\
&\approx \dfrac{\Delta x}{2}\sum_{i=1}^{N}\big(P^{n}(x_{i+1})+P^{n}(x_{i})\big) \leq 1 \qquad \forall n
\end{align*}

\noindent and thus the escape rate (per unit time), $r^{n}$, is

\begin{equation}
r^{n} = \dfrac{1 - s^{n}}{\Delta t}
\label{sec:SN escape}
\end{equation}

\noindent We will now consider two examples of bifurcation-induced tipping, one the classical normal form for the saddle-node and the other a model used to simulate the Indian monsoon. In these examples we will add noise to demonstrate the presence and behaviour of the early-warning indicators: increased variance and increased autocorrelation along with the escape rate all calculated using the formulas \eqref{sec:SN variance}, \eqref{autocorrelation SN sec} and \eqref{sec:SN escape} respectively. Linearising the bottom of the potential well $U(x)$ we can apply an approximation to the autocorrelation and variance (detailed later in Section \ref{subsec: lin drift}). We will also compare the numerical escape rate with Kramers' escape rate, $r_{K}$:

\begin{equation*}
r_{K} = \dfrac{\sqrt{\alpha\beta}}{2\pi}\exp\bigg(-\dfrac{\Delta U}{D}\bigg)
\end{equation*}

\noindent where $\alpha$ and $\beta$ represent the modulus of the curvature of the well and hill top of the potential and $\Delta U$ is the height of the potential barrier. Our first example is to consider slow passage towards a saddle-node bifurcation with two different types of drift.

\section{Slow passage towards a saddle-node bifurcation}
\label{sec: SN examples}

A saddle-node bifurcation can arise in 2 possible scenarios either a single fixed point appears and then splits into two fixed points that move further away from each other. The other scenario, which we will consider, is to start with two fixed points that move together,  eventually colliding and then disappear \citep{arrowsmith1992dynamical}. The following theorem from \citet{glendinning1994stability} provides further properties of the saddle-node bifurcation.

\begin{theorem}[Saddle-node bifurcation]

Suppose that for an ODE

\begin{equation*}
\dot{x} = G(x,p)
\end{equation*}

\noindent with $G(x_{*},p_{*}) = G_{x}(x_{*},p_{*}) = 0$. Then provided 

\begin{equation*}
G_{p}(x_{*},p_{*}) \neq 0 \qquad \mbox{and} \qquad G_{xx}(x_{*},p_{*}) \neq 0
\end{equation*}

\noindent there is a continuous curve of stationary points in a neighbourhood of $(x,p) = (x_{*},p_{*})$ which is tangent to the line $\{(x,p): p = p_{*}\}$ at the saddle-node bifurcation $(x_{*},p_{*})$. If $G_{p}G_{xx} < 0$ when evaluated at $(x_{*},p_{*})$ in some sufficiently small neighbourhood of $p = p_{*}$ then:

\begin{itemize}
\item there are no stationary points near $(x_{*},p_{*})$ if $p < p_{*}$
\item there are two stationary points near $x = x_{*}$ for each value of $p > p_{*}$
\item for $p \neq p_{*}$ both stationary points are hyperbolic and the lower branch is stable (unstable) and the upper branch is unstable (stable) if $G_{xx} > 0$ (or $G_{xx} < 0$) when evaluated at $(x_{*},p_{*})$
\end{itemize}

\noindent The first two statements are reversed if $G_{p}G_{xx} > 0$ when evaluated at $(x_{*},p_{*})$.
\label{SN theorem}

\end{theorem} 

All systems whose dynamics on the center manifold that can be described by Theorem \ref{SN theorem} at the equilibrium point $(x_{*},p_{*})$ are locally topologically equivalent to

\begin{equation}
\dot{x} = x^{2} - p
\label{SN normal form}
\end{equation}

\begin{figure}[h!]
        \centering
        \subcaptionbox{\label{SN bif}}[0.45\linewidth]
                {\includegraphics[scale = 0.3]{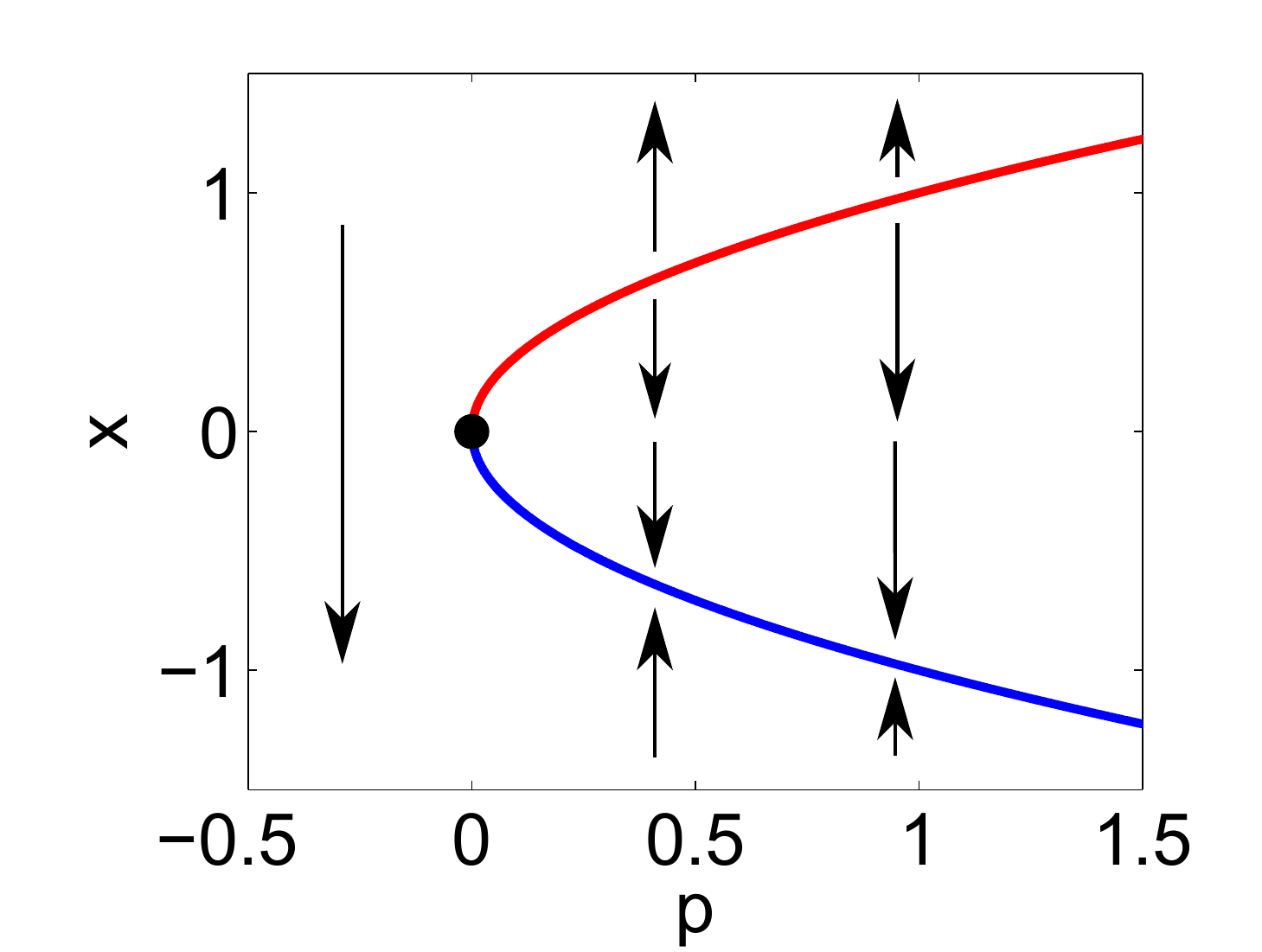}}
       \hfill
        \subcaptionbox{\label{SN bif2}}[0.45\linewidth]
                {\includegraphics[scale = 0.3]{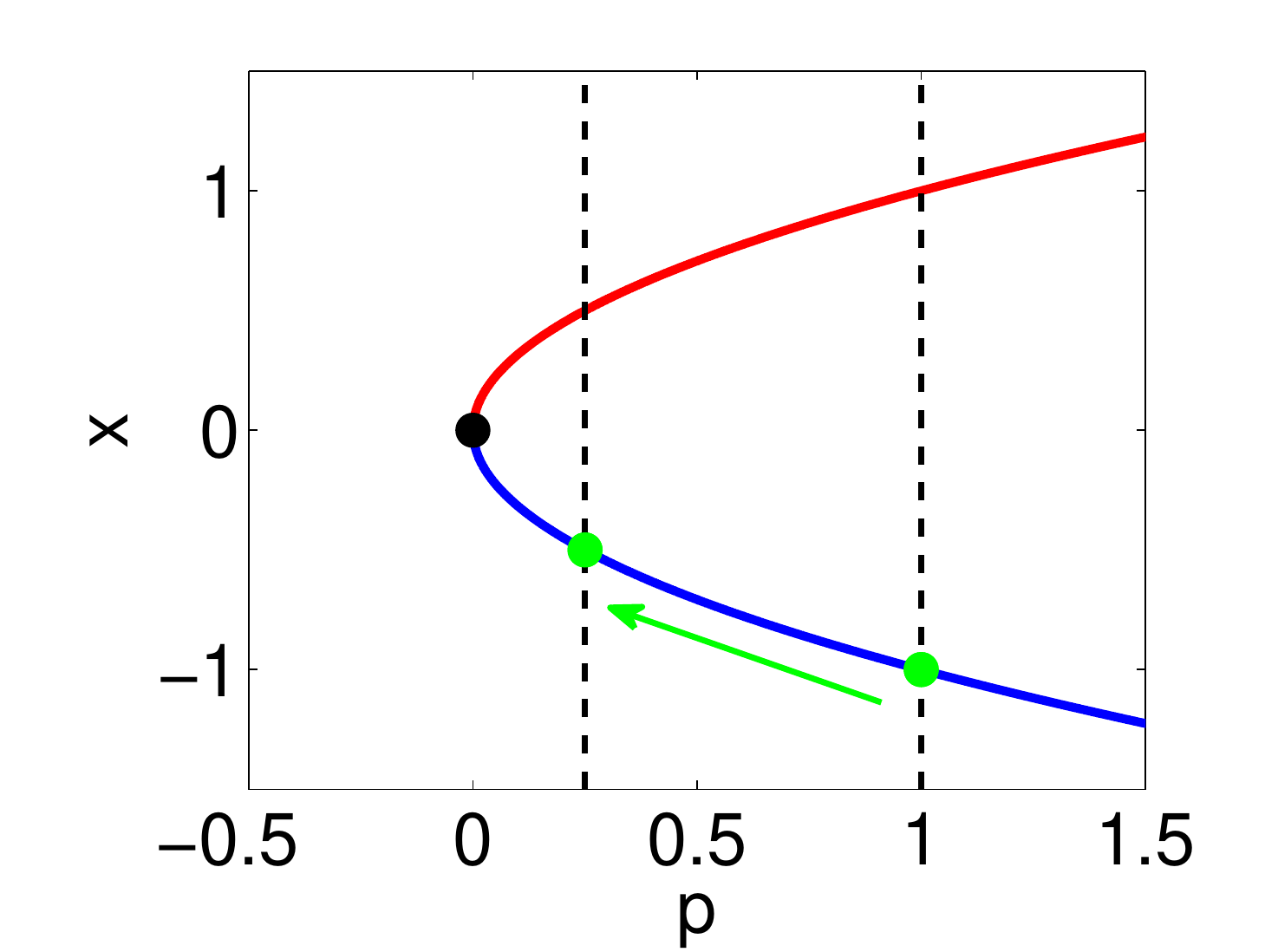}}
        ~ 
        \caption[Bifurcation diagram for saddle-node normal form.]{Bifurcation diagram for saddle-node normal form \eqref{SN normal form}. (a) Saddle-node bifurcation point at $(0,0)$, lower branch is stable (blue), upper branch unstable (red). (b) Illustration of slow passage towards bifurcation, $p(t) = p_{0} - \epsilon t$ starting at $p_{0} = 1$ and finishing at $p(T_{\mathrm{end}}) = 0.25$ where $\epsilon = 0.0075$ and $T_{\mathrm{end}} = 100$.}\label{sec:SNE Deterministic plots}
\end{figure}

\noindent \change{the normal (simplest) form for the saddle-node bifurcation. We will use equation \eqref{SN normal form} as our first example,} which has a saddle-node bifurcation at $(x,p) = (0,0)$. For \eqref{SN normal form} there exist no stationary points for $p < 0$ and two branch out from the bifurcation point for $p > 0$ see Figure \ref{SN bif}. The lower branch corresponding to $x_{s} = -\sqrt(p)$  is stable depicted in blue, and the upper branch, $x_{u} = +\sqrt{p}$ in red is unstable.  

\subsection{Linear drift}
\label{subsec: lin drift}

\begin{figure*}
\includegraphics[scale=0.35]{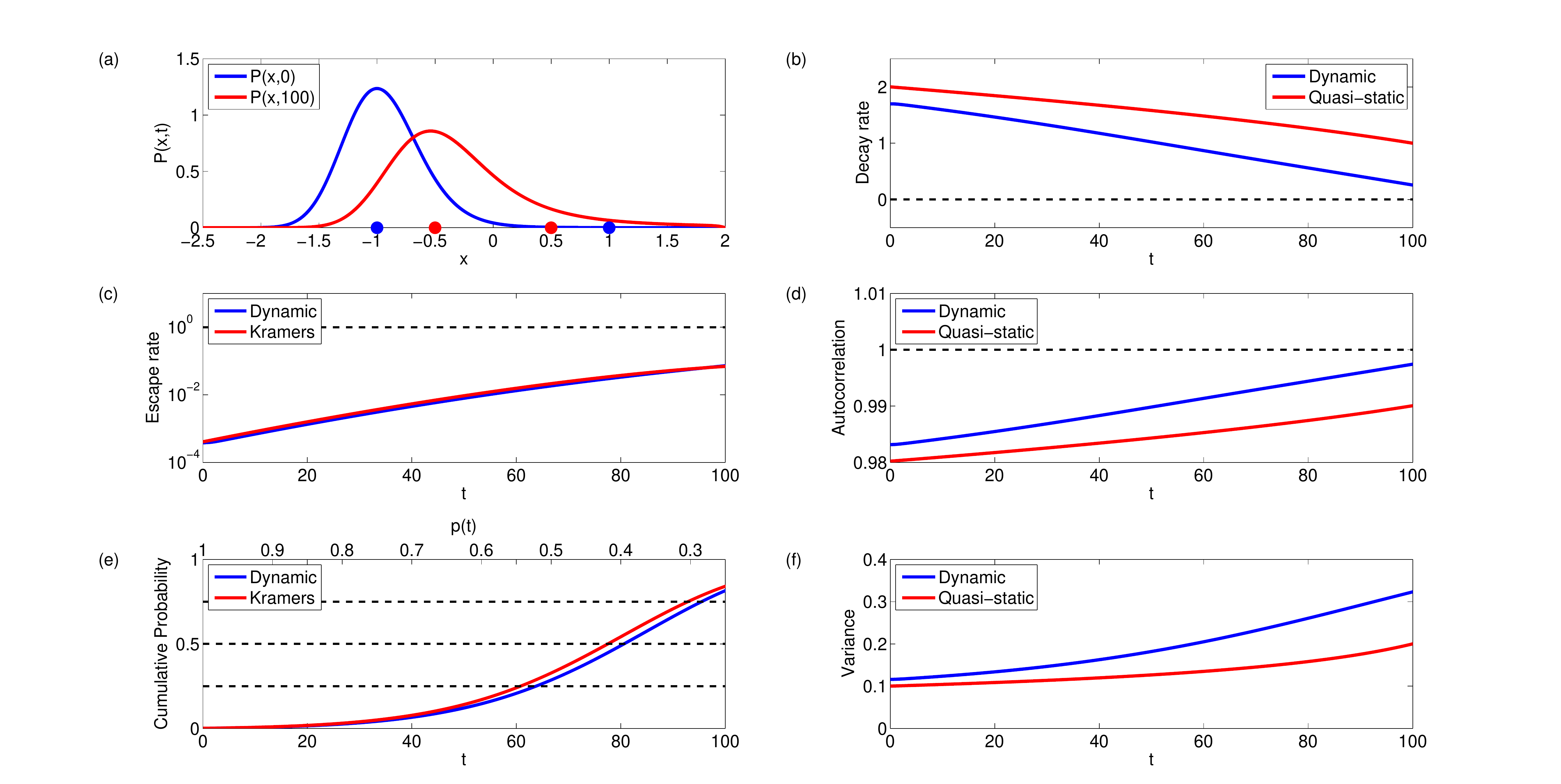}
\caption[Analysis of slow passage towards a saddle-node bifurcation with slow drift.]{Results of slow passage towards a saddle-node bifurcation according to SDE \eqref{SN SDE} and the corresponding Fokker-Planck equation \eqref{Fokker Planck}. Linear drift example: $p(t) = p_{0} - \epsilon t$, $p_{0} = 1$, $\epsilon = 0.0075$, $D = 0.2$, $\Delta x = 0.05$, $\Delta t = 0.01$. (a) Stationary probability density $P_{s}(x) = P(x,0)$ (blue) and final density $P(x,T_{\mathrm{end}})$ (red), fixed points indicated by dots on $x$-axis. Decay rate (b), lag-1 autocorrelation (per unit time) (d) and variance (f) for both nonlinear dynamic (numerical calculation, blue) and linear quasi-static (linearised Ornstein-Uhlenbeck, red). Escape rate (c) and Cumulative probability of escape (e) for dynamic (blue) and Kramers' (red).}
\label{SN eg1b labels}
\end{figure*}

We vary the bifurcation parameter, $p$, linearly according to

\begin{equation}
p \equiv p(t) = p_{0} - \epsilon t
\label{Linear drift}
\end{equation}

\noindent where $p_{0} = 1$. We fix the speed, $\epsilon = 0.0075$, at which we approach the saddle-node, which corresponds to a slow drift. We continue to $T_{\mathrm{end}} = 100$, which corresponds to $p(T_{\mathrm{end}}) = 0.25$, see Figure\ref{SN bif2}. Note that in this example we only approach the bifurcation and do not pass through it. This alteration is designed to allow us to detect the early-warning signals for an approach to a bifurcation-induced tipping event when white noise is added to \eqref{SN normal form}. \change{The choice of the parameter $p_{0} = 1$ means that for $t \approx 0$ we start in a slowly drifting well that is deep compared to the noise level.} The SDE for the random variable $X_{t}$ is:

\begin{equation}
\mathrm{d}X_{t} = \big(X_{t}^{2} - p(t)\big)\mathrm{d}t + \sqrt{2D}\mathrm{d}W_{t}
\label{SN SDE}
\end{equation}

\noindent and so by using the corresponding Fokker-Planck equation we can determine numerically the probability density function $P(x,t)$, as well as determine the early-warning indicators via the methods presented in Section \ref{sec: Numerical FP}. The results of this study are presented in Figure \ref{SN eg1b labels}.

Panel (a) displays the normalised stationary probability density $P_{s}(x) = P(x,0)$ in blue and the normalised probability density at $t = T_{\mathrm{end}}$ in red, on the fixed domain $[x_{\mathrm{start}},x_{\mathrm{end}}] = [-2.5,2]$. The stationary density has the form of a normal distribution that is centred about the stable point $x = -1$ (indicated by left blue dot), and the width (or variance) is determined by the noise level $D$. The final density in red is still centred close to the stable point $x_{s} = -\sqrt{p(T_{\mathrm{end}})}$ (left red dot) in the autonomous system \eqref{SN normal form}, because the drift $\epsilon$ is slow. However, the density has widened and gained a larger tail due to the shallowing of the potential well as the saddle-node bifurcation is approached.

Panels (c) and (e) present the escape rate and the cumulative probability of escape respectively. We define escape as a realisation having reached the upper boundary $x_{\mathrm{end}} = 2$. On the other hand, the lower boundary $x_{\mathrm{start}} = -2.5$ is set sufficiently far such that a realisation is extremely unlikely to reach this boundary. In both panels we compare the dynamic escape calculated numerically in blue with Kramers' escape rate in red. We can see that there is a good match between the two, initially there is a comparatively small escape but increases exponentially as we slowly move closer to the bifurcation point. Although when we compare for the cumulative probability of escape, panel (e), this demonstrates that Kramers' method slightly overestimates the escape compared to the numerical calculations. Moreover, we observe that even though the saddle-node bifurcation is not reached there is about an $80\%$ chance of incurring a bifurcation-induced with noise tipping event, where tipping refers to a realisation reaching $x_{\mathrm{end}}$. Note this is not a purely noise-induced tipping event as this would require the system to be stationary. Whereas, approaching the saddle-node weakens the stability of the steady state and therefore increases the vulnerability of the system tipping due to noise. We will now now consider the early-warning indicators to see if it is possible to detect the approach of the tipping point.

The numerically calculated lag-1 autocorrelation (per unit time) and variance of the probability densities (see Section \ref{subsec: Num EWI} for methods of calculation) are shown in blue in panels (d) and (f) respectively. We discover that both the autocorrelation and variance increase as the bifurcation is approached and therefore detects the possibility of a tipping event. In addition to the numerical calculations we can gain an approximation for the linear quasi-static autocorrelation and variance if we linearise around the potential well. As in the derivation of Kramers' escape rate we linearise by expanding $U(x)$ at the point $x_{s}$ and so we get

\begin{equation*}
U(x) =  U(x_{s}) + \dfrac{1}{2}\kappa(x - x_{s})^{2} + \mathcal{O}(x-x_{s})^{3}
\end{equation*}

\noindent where $\kappa = U''(x_{s})$ corresponds to the decay rate. The linearised SDE is the Ornstein-Uhlenbeck process for the random variable $X_{t}$:

\begin{equation}
\mathrm{d}X_{t} = -\kappa X_{t}\mathrm{d}t + \sqrt{2D}\mathrm{d}W_{t}
\label{sec:SNE OU SDE}
\end{equation} 

\noindent For the static Ornstein-Uhlenbeck process the autocorrelation and variance is given by \citep{aalen2008survival}:

\begin{align}
\nonumber
\text{Autocorrelation:} \qquad &a = \exp(-\kappa\Delta t) \approx 1 - \kappa\Delta t \\
\label{OU variance}
\text{Variance:} \qquad &V = \dfrac{D}{\kappa}
\end{align}

\noindent Therefore, at the end of each time step $\Delta t$ we calculate the linear quasi-static decay rate (b), autocorrelation (d) and variance (f) in red. Once again the indicators, autocorrelation and variance increase, as we move closer to the bifurcation, indicating a possible bifurcation-induced tipping event. However, there is a clear difference between the linear quasi-static and the nonlinear dynamic calculated values. \change{The following section analyses this systematic difference in detail.}

\subsubsection*{Systematic differences between nonlinear dynamic and linear quasi-static values of decay rate and early-warning indicators}

The notable differences between the true (nonlinear dynamic) and approximate (linear quasi-static) values of the decay rate and early-warning indicators are:

\begin{itemize}
\item \textbf{Higher order terms in potential:} we apply a linearisation about the stable equilibrium in the autonomous system \eqref{SN normal form} to use the Ornstein-Uhlenbeck process (Spatial error in $x$)
\item \textbf{Time dependence of process:} Ornstein-Uhlenbeck is a stationary process whereas we consider a slow passage towards a saddle-node bifurcation. (Time error in $t$)
\end{itemize}

We again state the SDE for the saddle-node bifurcation \eqref{SN SDE} coupled with the ODE for the bifurcation parameter for which the solution \eqref{Linear drift} satisfies with initial condition $p(0) = p_{0}$:

\begin{align}
\label{x SDE}
\mathrm{d}X_{t} &= (X_{t}^{2} - p)\mathrm{d}t + \sqrt{2D}\mathrm{d}W_{t} \\
\mathrm{d}p &= -\epsilon\mathrm{d}t
\label{p ODE}
\end{align}

\noindent Scaling system \eqref{x SDE}--\eqref{p ODE} will allow us to establish the parameters we need to vary to analyse the impacts of the two points mentioned above. We apply the following scalings:

\begin{equation*}
X = \alpha Y, \qquad t = \beta\tau, \qquad p = \gamma q
\end{equation*}

\noindent substituting these scalings into \eqref{x SDE}--\eqref{p ODE} gives:

\begin{align*}
\alpha\mathrm{d}Y_{\tau} &= (\alpha^{2}Y_{\tau}^{2} - \gamma q)\beta\mathrm{d}\tau + \sqrt{2D\beta}\mathrm{d}W_{\tau} \\
\gamma\mathrm{d}q &= -\epsilon\beta\tau
\end{align*}

\noindent and rearranging produces

\begin{align}
\label{scaled SDE}
\mathrm{d}Y_{\tau} &= \bigg(\alpha\beta Y_{\tau}^{2} - \dfrac{\gamma\beta}{\alpha}q\bigg)\mathrm{d}\tau + \sqrt{\dfrac{2D\beta}{\alpha^{2}}}\mathrm{d}W_{\tau} \\
\mathrm{d}q &= -\dfrac{\epsilon\beta}{\gamma}\mathrm{d}\tau
\label{scaled ODE}
\end{align}

\noindent Setting 

\begin{equation*}
\alpha\beta = 1, \qquad \dfrac{\gamma\beta}{\alpha} = 1, \qquad \dfrac{D\beta}{\alpha^{2}} = \tilde{D}, \qquad \tilde{\epsilon} = \dfrac{\epsilon\tilde{D}}{D}
\end{equation*}

\noindent reduces the scaled system \eqref{scaled SDE}--\eqref{scaled ODE} to be of the form of \eqref{x SDE}--\eqref{p ODE}:

\begin{align}
\label{y SDE}
\mathrm{d}Y_{\tau} &= (Y_{\tau}^{2} - q)\mathrm{d}\tau + \sqrt{2\tilde{D}}\mathrm{d}W_{\tau} \\
\mathrm{d}q &= -\tilde{\epsilon}\mathrm{d}\tau
\label{q ODE}
\end{align}

\noindent where the scalings are given as

\begin{equation}
X = \bigg(\dfrac{D}{\tilde{D}}\bigg)^{\frac{1}{3}}Y, \quad t = \bigg(\dfrac{\tilde{D}}{D}\bigg)^{\frac{1}{3}}\tau, \quad p = \bigg(\dfrac{D}{\tilde{D}}\bigg)^{\frac{2}{3}}q
\label{scalings}
\end{equation}

\noindent Following the rescaling we are in a position such that we can fix one of the following parameters: the noise level $D$, drift speed $\epsilon$ or $p_{0}$ which is proportional to the distance from the saddle-node.

\paragraph{Higher order terms in potential:} We will initially consider zero drift speed ($\epsilon = 0$) and analyse the impact the linearisation has by considering the nonlinear and linear decay rate and variance for different noise levels $D$. \change{We choose the decay rate instead of autocorrelation because the decay rate is independent of the time step.} We set $p(t) = p_{0} = 1$ in \eqref{SN SDE} without loss of generality (w.l.o.g.) as a rescaling of space and time demonstrated that considering a $q_{0}$ different from $p_{0}$ would equate to using a different noise level $\tilde{D}$ given by the relation:

\begin{equation}
\tilde{D} = \bigg(\dfrac{q_{0}}{p_{0}}\bigg)^{\frac{3}{2}}D
\label{noise scaling}
\end{equation}

\begin{figure}[h!]
        \centering
        \subcaptionbox{\label{SN decay D}}[0.45\linewidth]
                {\includegraphics[scale = 0.3]{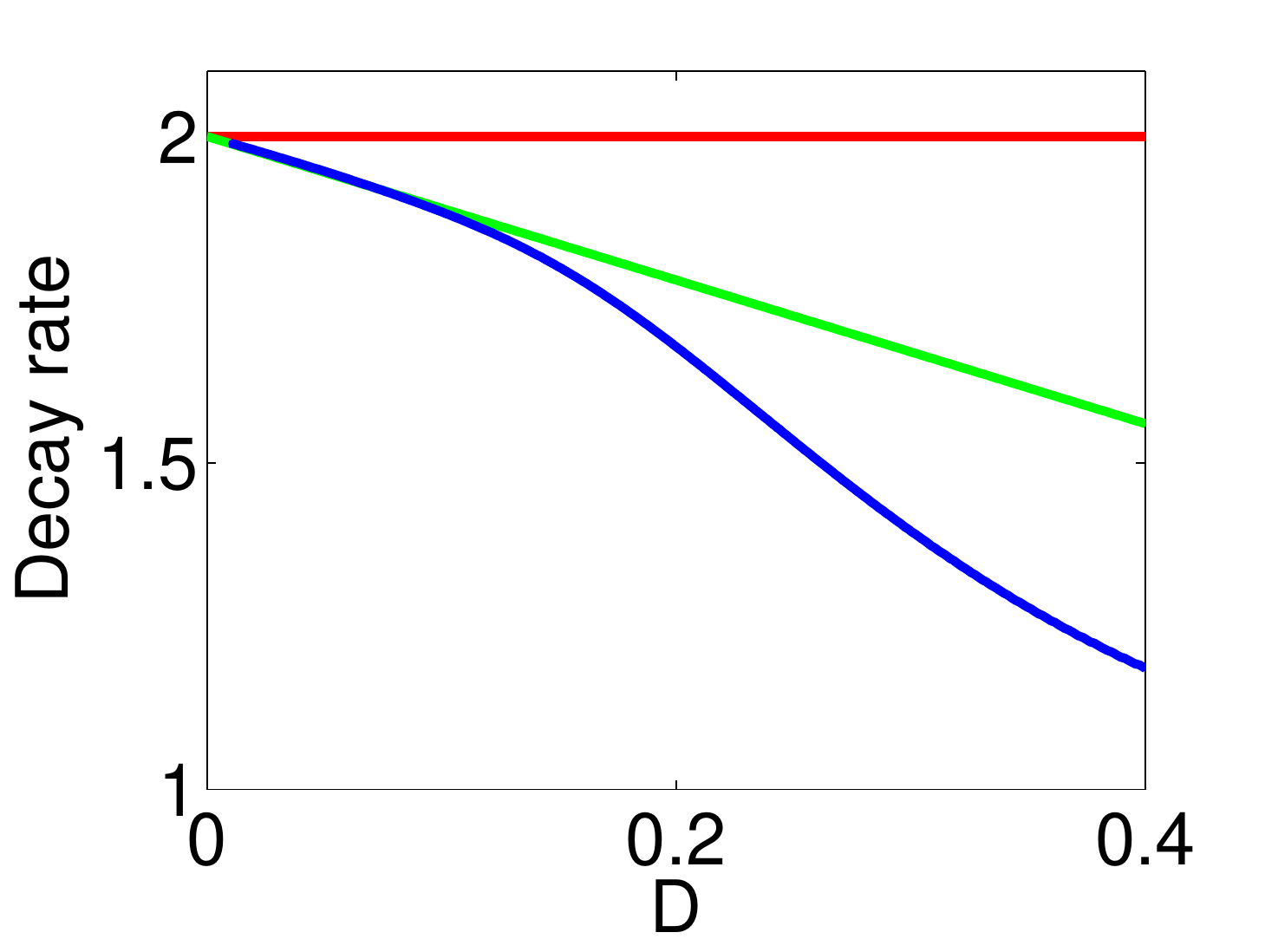}}
       \hfill 
        \subcaptionbox{\label{SN variance D}}[0.45\linewidth]
                {\includegraphics[scale = 0.3]{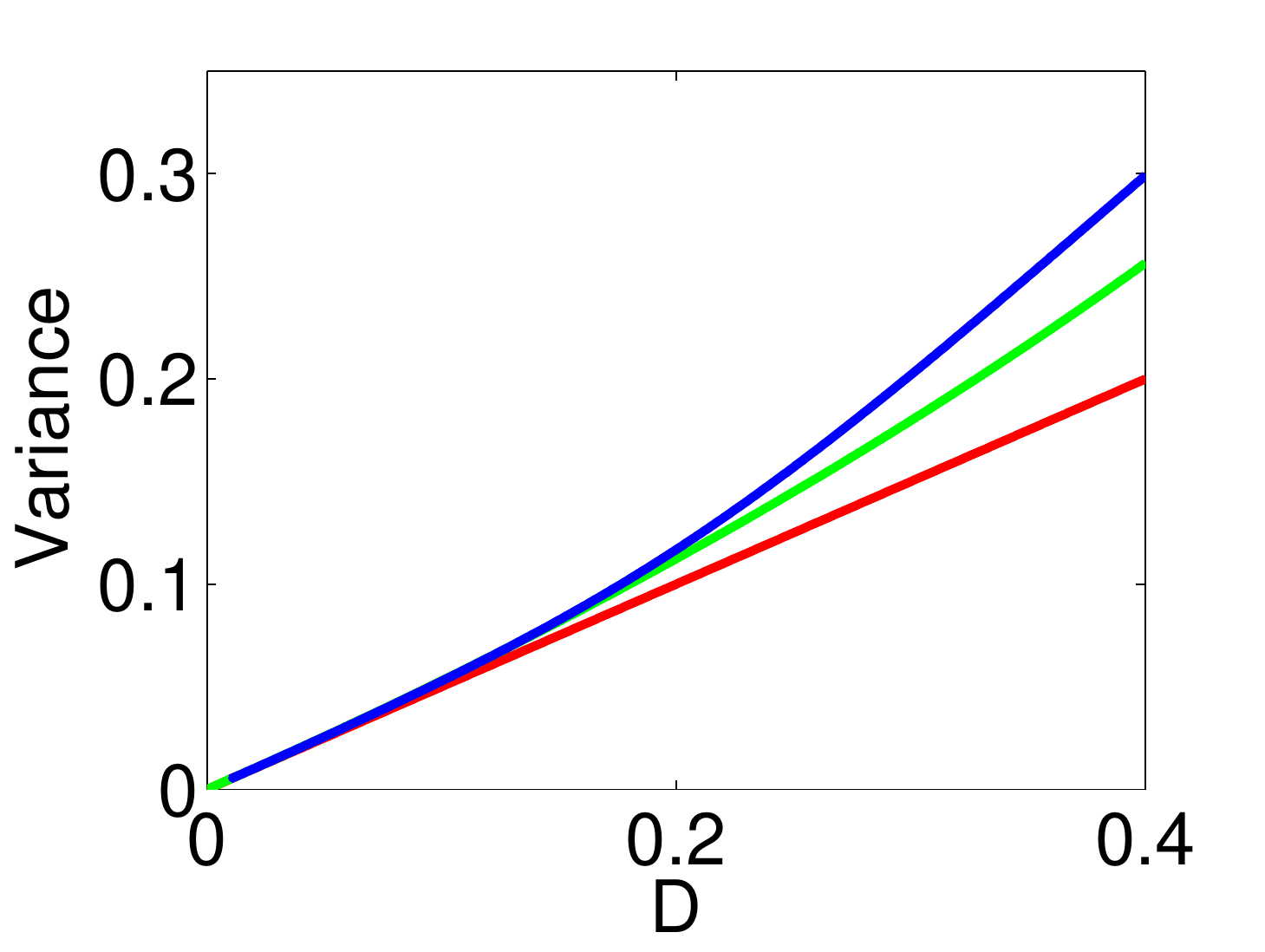}}
        \caption[Comparison of nonlinear and linear decay rate and variance.]{Comparison of nonlinear (blue) decay rate (a) and variance (b) emanating from the SDE \eqref{SN SDE} with the linear versions (red) \eqref{sec:SNE OU SDE} for a range of noise levels $D$ with linear drift $p(t) = p_{0} = 1$ $(\epsilon = 0)$. Green curve in (a), $\kappa_{n}$, \eqref{approx dynamic DR} linear fit to the nonlinear decay rate for small $D$ and corresponding variance (b) $V_{n}$, given by \eqref{approx dynamic Var}.}\label{Q vs C}
\end{figure}

In Figure \ref{SN decay D}, we present the nonlinear (blue) and linear (red) decay rates for a range of noise levels. The linear decay rate $\kappa$ is independent of $D$ and thus remains a constant value. For small noise $D<0.1$, the nonlinear decay rate $\kappa_{n}$ follows a linear trend, roughly approximated by the formula:

\begin{equation}
\kappa_{n} = \kappa_{l} - \dfrac{\kappa_{c}\tilde{D}}{q_{0}}
\label{approx dynamic DR}
\end{equation}

\noindent where $\kappa_{l}$ is the linear decay rate $\kappa$ from the Ornstein-Uhlenbeck process \eqref{sec:SNE OU SDE} and $\kappa_{c} \approx 1.1$, a constant determined by the tangent to the initial nonlinear decay rate. This tangent is plotted in green in Figure \ref{SN decay D} and provides a good fit for $D \leq 0.1$. Though the linear fit loses accuracy for larger noise levels as the nonlinear decay rate decreases quicker. 

The nonlinear (blue) and linear (red) variances for the range of noise levels $D$ is displayed in Figure \ref{SN variance D}. The linear variance gives a linear trend since the linear decay rate is constant, see equation \eqref{OU variance}. The linear variance underestimates the nonlinear variance, especially for large noise. The green curve represents the linear fit of the nonlinear decay rate \eqref{approx dynamic DR} applied to the linear formula \eqref{OU variance} for the variance:

\begin{equation}
V_{n} = \dfrac{D}{\kappa_{n}}
\label{approx dynamic Var}
\end{equation}  

\noindent This gives a much improved approximation to the nonlinear variance for small noise.     

\paragraph{Time dependence of process:} For zero drift speed, we fixed $p_{0}$ w.l.o.g. and varied the noise level $D$. One reason for this, was to allow us to easily fit a linear trend to the nonlinear decay rate for small $D$. Considering a non-zero drift speed it is advantageous to set $D = 0.2$ w.l.o.g. and leave $p_{0}$ and $\epsilon$ free. \change{We choose $p_{0}$ large, such that the error between the nonlinear and linear decay rate and variance is minimal, for example $p_{0} = 4$, which is equivalent to $q_{0} = 1$, $D = 0.025$, }

\begin{figure}[h!]
        \centering
          \subcaptionbox{\label{SN decay drift small D}}[0.45\linewidth]
                {\includegraphics[scale = 0.3]{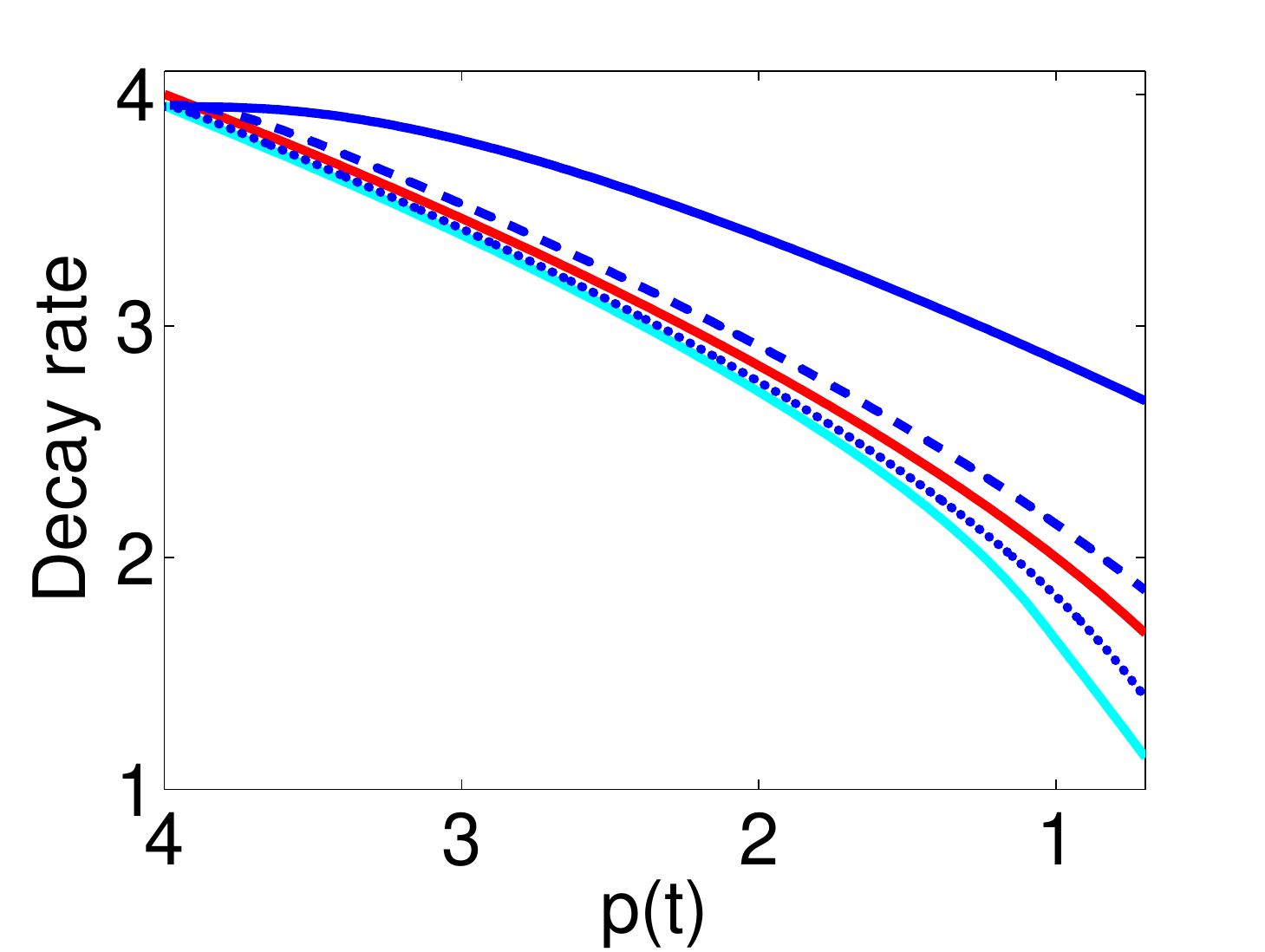}}
       \hfill
        \subcaptionbox{\label{SN variance drift small D}}[0.45\linewidth]
                {\includegraphics[scale = 0.3]{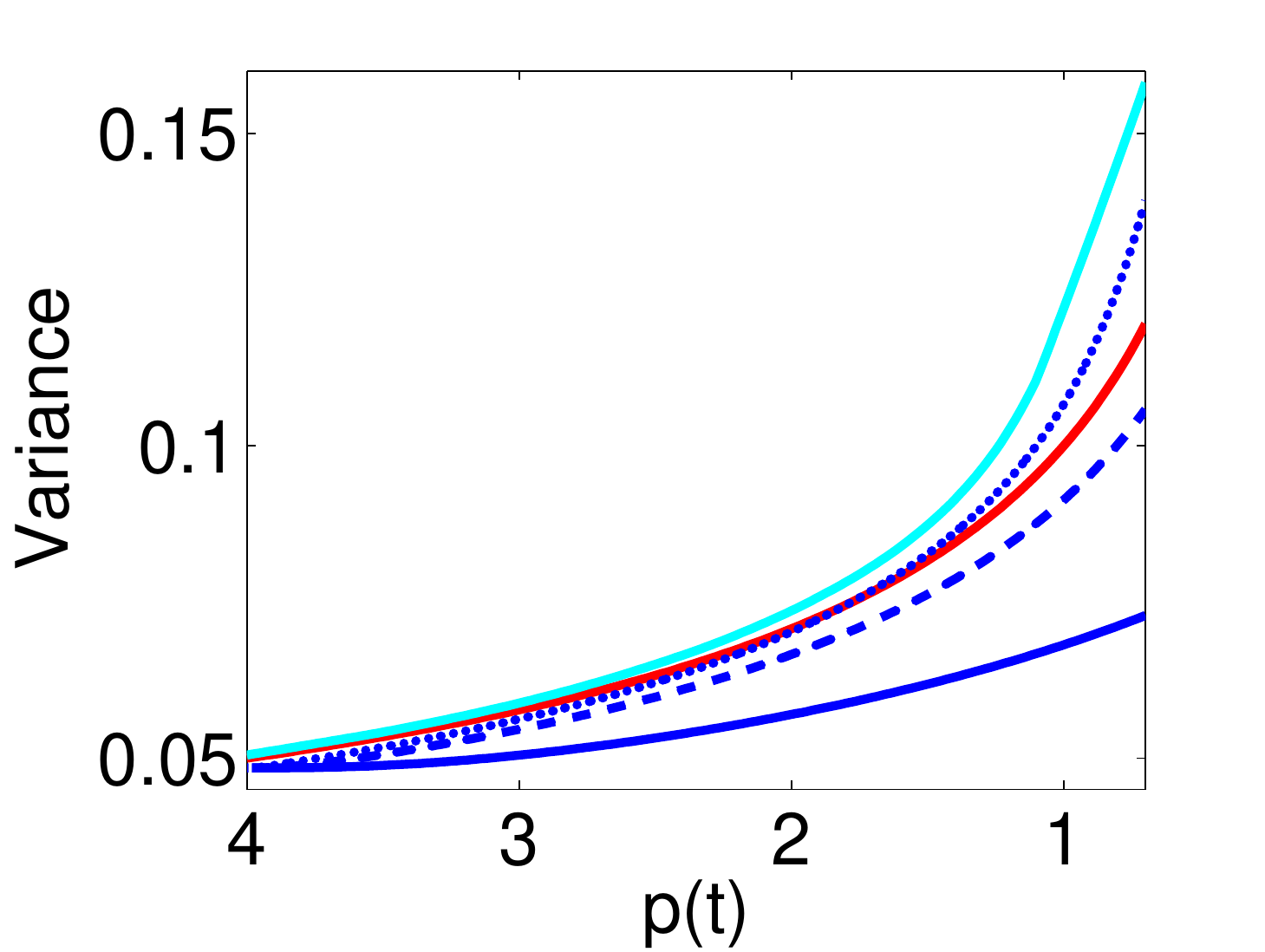}}
        \caption[Drift speed analysis on nonlinear dynamic decay rate and variance compared with linear and nonlinear quasi-static approximations.]{Assessing the impact different drift speeds $\epsilon$ has on the decay rate (a) and variance (b), for noise level $D = 0.2$. Nonlinear dynamic decay rate and variance for \eqref{x SDE}--\eqref{p ODE} given in dark blue for different drift speeds; $\epsilon = 2.5$ (solid), $\epsilon = 0.5$ (dashed), $\epsilon = 0.1$ (dotted). Linear quasi-static approximations for \eqref{sec:SNE OU SDE} given in red. Light blue curves represent the nonlinear quasi-static approximations \eqref{nonlin qstatic decay} and \eqref{nonlin qstatic var}, used to indicate the error caused solely by the drift speed $\epsilon$.}\label{Drift analysis}
\end{figure}

Figure \ref{Drift analysis} \change{provides a comparison between the nonlinear dynamic, nonlinear quasi-static and linear quasi-static approximations for different drift speeds $\epsilon$ (time lengths adjusted accordingly).} The fastest drift speed $\epsilon = 2.5$ is given in solid blue, $\epsilon = 0.5$ in dashed blue and the slowest drift $\epsilon = 0.1$ is given by dotted blue. The red curve gives the linear quasi-static approximation for each fixed time point. \change{The light blue curve is the nonlinear quasi-static approximation, designed to indicate the error caused solely by the drift speed $\epsilon$. Let us briefly describe how the nonlinear quasi-static approximation is calculated.} 

\change{If we treat $t$ as a parameter then for every fixed $p(t)$ and noise level $D$ in system \eqref{x SDE}--\eqref{p ODE}, is equivalent to considering in the scaled system \eqref{y SDE}--\eqref{q ODE} ($\tilde{\epsilon} = 0$), $q_{0} = 1$ with noise level $\tilde{D}$ given by the relation \eqref{noise scaling}.} For example, $p(t_{0}) = p_{0} = 4$ and $D = 0.2$ this is equivalent to considering $q_{0} = 1$ with a noise level $\tilde{D} = 0.025$. \change{For $q_{0} = 1$ we have previously calculated the nonlinear decay rate and variance for a range of noise levels $\tilde{D}$, see Figure \ref{Q vs C}. Though, we have not yet established the relationship for the decay rate and variance between the two systems. We use the linearised Ornstein-Uhlenbeck process to determine the scalings of the decay rate and variance. For the original system \eqref{x SDE}--\eqref{p ODE} has decay rate $\kappa_{X} = 2\sqrt{p}$ and the scaled system of equations \eqref{y SDE}--\eqref{q ODE} has decay rate $\kappa_{Y} = 2\sqrt{q}$}. We can therefore express the nonlinear decay rate for \eqref{x SDE}--\eqref{p ODE}, $\kappa_{X}$ in terms of the nonlinear decay rate for the scaled system \eqref{y SDE}--\eqref{q ODE}, $\kappa_{Y}$, which we know, using the scaling between $p$ and $q$ given in \eqref{scalings}:

\begin{equation}
\kappa_{X} = 2\sqrt{p} = 2\sqrt{q}\bigg(\dfrac{D}{\tilde{D}}\bigg)^{\frac{1}{3}} = \bigg(\dfrac{D}{\tilde{D}}\bigg)^{\frac{1}{3}}\kappa_{Y}
\label{nonlin qstatic decay}
\end{equation}

\noindent Likewise, for the variance:

\begin{equation}
V_{X} = \dfrac{D}{\kappa_{X}} = \dfrac{D}{\kappa_{Y}}\bigg(\dfrac{\tilde{D}}{D}\bigg)^{\frac{1}{3}} = \dfrac{\tilde{D}}{\kappa_{Y}}\bigg(\dfrac{D}{\tilde{D}}\bigg)^{\frac{2}{3}} = \bigg(\dfrac{D}{\tilde{D}}\bigg)^{\frac{2}{3}}V_{Y}
\label{nonlin qstatic var}
\end{equation}

\noindent Applying these scalings will produce an approximation for the nonlinear quasi-static decay rate and variance that takes into consideration higher order terms (softening) of the potential but not the drift speed.

Figure \ref{Drift analysis} presents the results of the drift analysis for the decay rate and variance respectively. For zero drift speed we observed little error between the nonlinear and linear decay rate and variance for small noise (c.f. Figure \ref{Q vs C}). Therefore the difference between the quasi-static nonlinear (light blue) and linear (red) approximations is small but increases slightly as the bifurcation is approached as this is equivalent to the noise range considered for constant drift when $p_{0} = 1$. For a fast drift, $\epsilon = 2.5$ both the decay rate and variance develop a large difference between the nonlinear dynamic (solid blue) and the quasi-static nonlinear or linear curves. Furthermore, we highlight the drift has the opposite effect linearisation had on the nonlinear decay rate and variance compared with the linear versions. The softening of the potential brought the nonlinear decay rate below the linear. However, the drift causes the nonlinear dynamic decay rate to be higher than the linear quasi-static. The opposite holds for the variance, i.e. the nonlinear is higher due to the softening and the drift brings the nonlinear dynamic variance below the linear quasi-static. For each curve of the nonlinear dynamic decay rate or variance we decrease the drift speed by a factor of $5$. Thus, the dashed blue corresponds to the nonlinear dynamic decay rate or variance for a drift speed, $\epsilon = 0.5$. The nonlinear dynamic decay rate and variance for this drift speed both start to converge towards the quasi-static curves. At the slowest drift speed considered, $\epsilon = 0.1$ (dotted blue) the nonlinear dynamic has moved below (above) the linear quasi-static for the decay rate (variance). This indicates that for $\epsilon = 0.1$ the error from linearisation dominates the drift speed error, whereas, for the previous two drift speeds the error from the drift dominated. 

We would like to emphasise that it is not possible to single out that either the linearisation or drift is the dominating factor for large noise and large drift speed for example. This is because the linearisation accounts for a spatial error, whereas the drift is a time error and it is therefore difficult to quantify what is `large' noise and equally `large' drift speed.    

\paragraph{Summary:} \change{We have discussed in detail two systematic differences between the nonlinear dynamic and linear quasi-static decay rates and variances, the results are summarised in Table \ref{Analysis summary}. One source of difference is linearisation of the potential well (nonlinear vs linear), this effect can be seen for zero drift ($\epsilon = 0$) and varying one of the parameters $D$ or $p_{0}$. The scaling \eqref{noise scaling} notifies us that increasing the noise level is equivalent to decreasing the parameter $p_{0}$. An increase in the noise level causes a greater overestimation ($+$) of the linear decay rate compared to the nonlinear decay rate. Consequently the linear autocorrelation underestimates ($-$) the nonlinear autocorrelation. The linear variance also provides an underestimation compared to nonlinear variance.}

\begin{table}
\begin{center}
\begin{ruledtabular}
    \begin{tabular}{cccc}
    Parameter & Decay & Autocorrelation & Variance\\[0.5ex] \hline
    {$\!\begin{aligned}\ \\[-3ex]   
               &D\uparrow \\
               &p_{0}\downarrow \end{aligned}$} &  {$\!\begin{aligned}\ \\[-3ex]   
               + \\
               + \end{aligned}$} &  {$\!\begin{aligned}\ \\[-3ex]   
               - \\
               - \end{aligned}$} &  {$\!\begin{aligned}\ \\[-3ex] 
               - \\
               - \end{aligned}$}  \\ \hline\ \\[-2ex]
               $\epsilon\uparrow$ & $-$ & $+$ & $+$\\ 
    \end{tabular}
    \caption[Summarising the effect linearisation and drift has on the decay rate and early-warning indicators.]{Summarising the effect linearisation and drift has on the decay rate and early-warning indicators autocorrelation and variance. For an increase of $D$ or decrease of $p_{0}$ a ($+$) implies the linear quantity overestimates the nonlinear quantity, whereas $(-)$ is an underestimation. For an increase of $\epsilon$ a ($+$) implies the quasi-static value overestimates the dynamic value, whereas $(-)$ is an underestimation.}
    \label{Analysis summary}
    \end{ruledtabular}
\end{center}
\end{table}

\change{We also investigated the time dependence of the process (dynamic vs quasi-static) by varying the drift speed $\epsilon$. For faster drift speeds $\epsilon$ we found that the nonlinear quasi-static decay rate underestimated the nonlinear dynamic decay rate. Whereas, for the early-warning indicators; autocorrelation and variance the nonlinear quasi-static overestimated the nonlinear dynamic indicators.}

\change{One should be aware it is difficult to quantify the error for considering both systematic differences together (nonlinear dynamic vs linear quasi-static). The linearisation of the potential well causes a spatial error, while, a time error is inflicted for non-zero drift speeds.}

\subsection{Nonlinear drift}
\label{subsec: SN nlin drift}

In the previous example, we used a linear drift which represented a slow passage towards the saddle-node bifurcation. For this example, we will consider a nonlinear drift that has been motivated by rate-induced tipping. 

Rate-induced tipping can be observed in the following ODE \citep{ashwin2012tipping}:

\begin{equation}
\dot{x} = (x + \lambda(t))^{2} - p_{0}
\label{sec:SNE rtip ODE}
\end{equation}

\noindent where

\begin{equation}
\lambda(t) = \dfrac{\lambda_{\max}}{2}\bigg(\tanh\bigg(\dfrac{\lambda_{\max}\epsilon t}{2}\bigg) + 1\bigg)
\label{sec:SNE lambda eq}
\end{equation}

\noindent is often referred to as a ramping parameter. Notice, that \eqref{sec:SNE rtip ODE} is the saddle-node normal form if we set $\lambda = 0$ and therefore as in the previous example $p_{0}$ determines the distance between the stable and unstable equilibria. The role of $\lambda(t)$ is to apply a shift to the saddle-node, governed by \eqref{sec:SNE lambda eq}, where $\lambda_{\max}$ determines the distance of the shift and $\epsilon$ is directly proportional to the rate of the shift. If we impose the condition on the shift parameter $\lambda(t)$, namely $\lambda(0) = \lambda_{\max}/2$ then \eqref{sec:SNE lambda eq} is a solution of the ODE:

\begin{equation}
\dot{\lambda} = \epsilon\lambda(\lambda_{\max} - \lambda)
\label{sec:SNE lambda ODE}
\end{equation} 

\noindent Applying a change of coordinates $y = x + \lambda$ to \eqref{sec:SNE rtip ODE} and using \eqref{sec:SNE lambda ODE} we have

\begin{align*}
\dot{y} &= \dot{x} + \dot{\lambda} \\
&= (x + \lambda)^{2} - p_{0} + \epsilon\lambda(\lambda_{\max} - \lambda) \\
&= y^{2} - (p_{0} - \epsilon\lambda_{\max}\lambda + \epsilon\lambda^{2})
\end{align*}  

\noindent which is simply the normal form of the saddle-node:

\begin{equation}
\dot{y} = y^{2} - \tilde{p}(t)
\label{transformed ODE}
\end{equation} 

\noindent where

\begin{equation}
\tilde{p}(t) = p_{0} - \epsilon\lambda_{\max}\lambda + \epsilon\lambda^{2}
\label{transformed drift}
\end{equation}

\noindent Though, this does not reduce rate-induced tipping to the previous case of bifurcation-induced tipping with noise. Recall, bifurcation-induced tipping with noise considers simply approaching or crossing a bifurcation transversally. Whereas, $\tilde{p}(t)$ represents a nonlinear drift, that causes the system to move towards the bifurcation but then retreats away from the bifurcation. Initially the system starts a distance $\sqrt{p_{0}}$ from the saddle-node bifurcation at $t = -\infty$. The speed at which the saddle-node is approached (or even past, depending on the choice of parameters) increases until $t = 0$. At $t = 0$ the system is the closest it gets to the saddle-node bifurcation (or furthest past the bifurcation). The drift then turns around and moves the system away from the bifurcation at a mirrored speed to which the approach occurred.

\begin{figure}[h!]
        \centering
          \subcaptionbox{\label{Transformed rtip drift}}[0.45\linewidth]
                {\includegraphics[scale = 0.3]{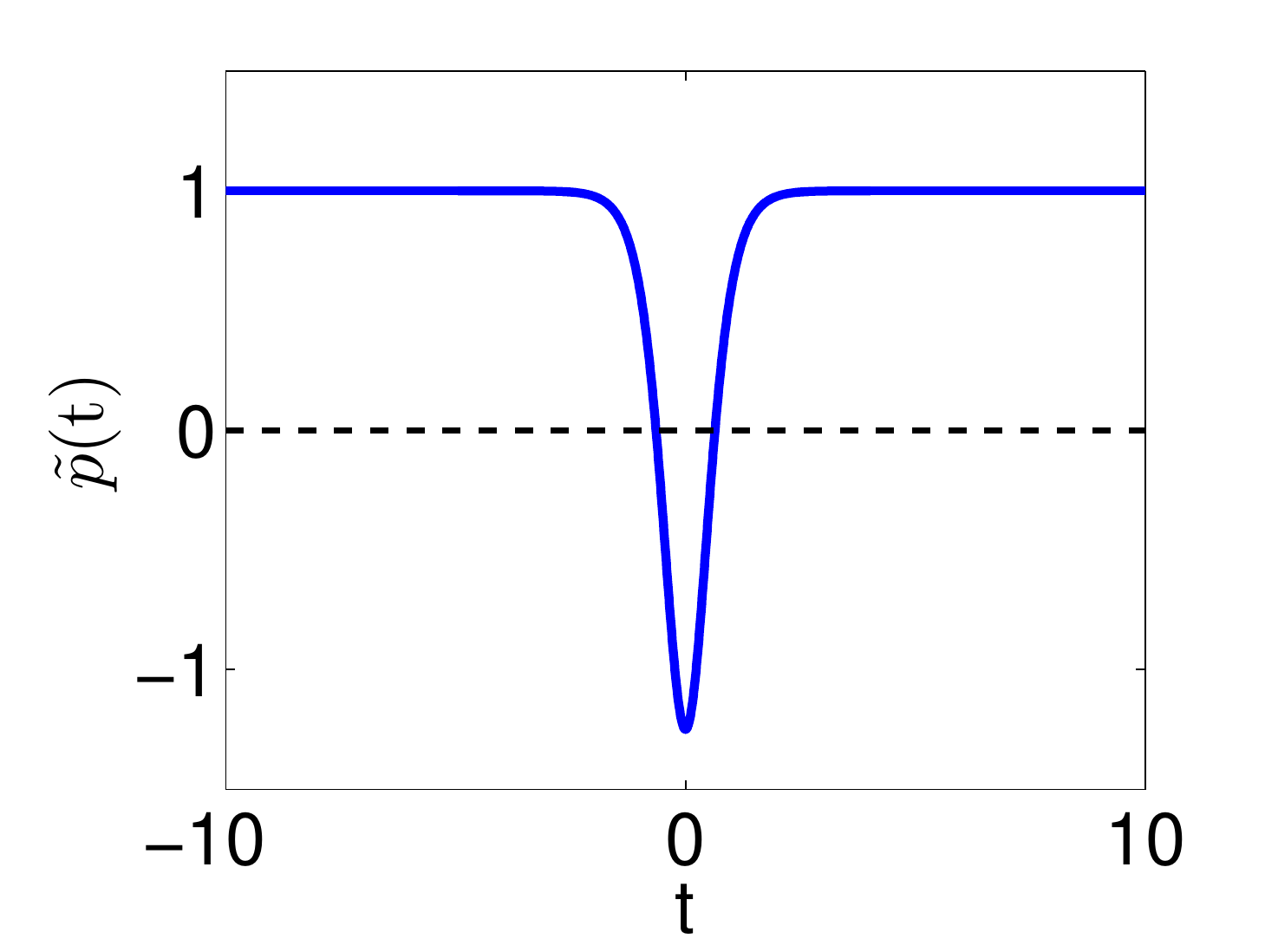}}
       \hfill
        \subcaptionbox{\label{Transformed rtip TP}}[0.45\linewidth]
                {\includegraphics[scale = 0.3]{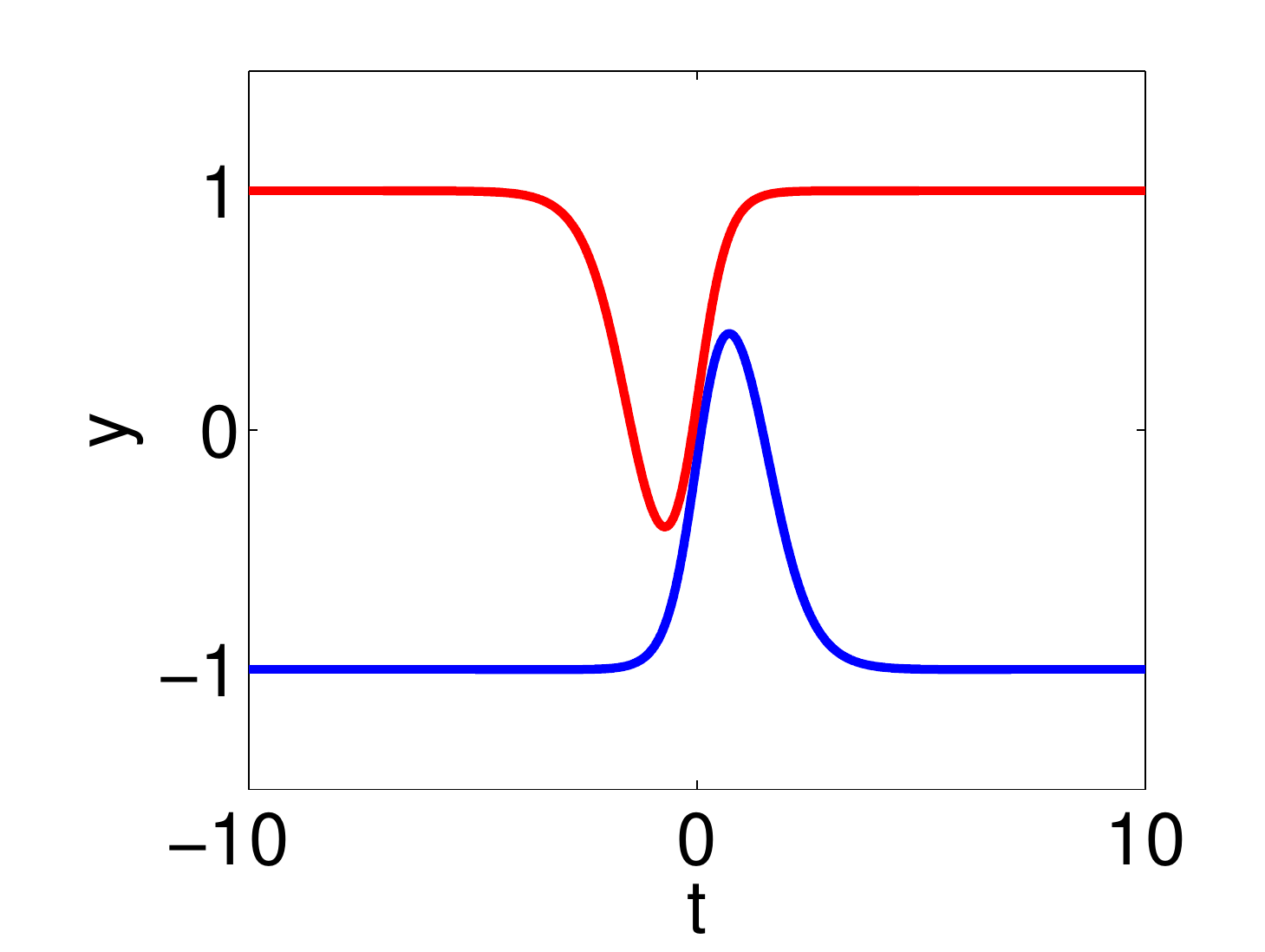}}
        \caption[Rate-induced system equivalent to saddle-node normal form with nonlinear drift by a change of coordinates.]{(a) Applying a change of coordinates to the rate-induced system \eqref{sec:SNE rtip ODE}, \eqref{sec:SNE lambda ODE} gives rise to a nonlinear drift $\tilde{p}(t)$ \eqref{transformed drift} plotted for the saddle-node normal form \eqref{transformed ODE}. Black dashed line represents saddle-node bifurcation. (b) Time profile for trajectories of a realisation that either starts at the stable quasi-static equilibrium (blue) or finishes at the unstable quasi-static equilibrium (red) for the system \eqref{transformed ODE}--\eqref{transformed drift}. Parameters: $p_{0} = 1$, $\lambda_{\max} = 3$, $\epsilon = 1$.}\label{Transformed rtip}
\end{figure}

Setting the parameters $\lambda_{\max} = 3$ and $p_{0} = 1$, rate-induced tipping occurs for drift speeds $\epsilon > 4/3$ \citep{perryman2015tipping}. Notice though, choosing $\epsilon = 1$ will result in passing the saddle-node bifurcation for a brief period of time, Figure \ref{Transformed rtip drift}. However, since $\epsilon = 1 < 4/3$, tipping does not occur and therefore the system recovers having passed the bifurcation. This is highlighted further in Figure \ref{Transformed rtip TP}, which displays the time profiles for a trajectory starting at $(x,t) = (-1,-10)$ (blue) and for one ending at $(x,t) = (1,10)$ (red). The red curve acts as a separatrix such that trajectories below will get attracted towards the blue curve and trajectories above will escape to $\infty$ in finite time. The blue curve demonstrates a transition past the saddle-node bifurcation before returning to the stable quasi-static equilibrium. 

\begin{figure}[h!]
        \centering
          \subcaptionbox{\label{SN nonlin drift}}[0.45\linewidth]
                {\includegraphics[scale = 0.3]{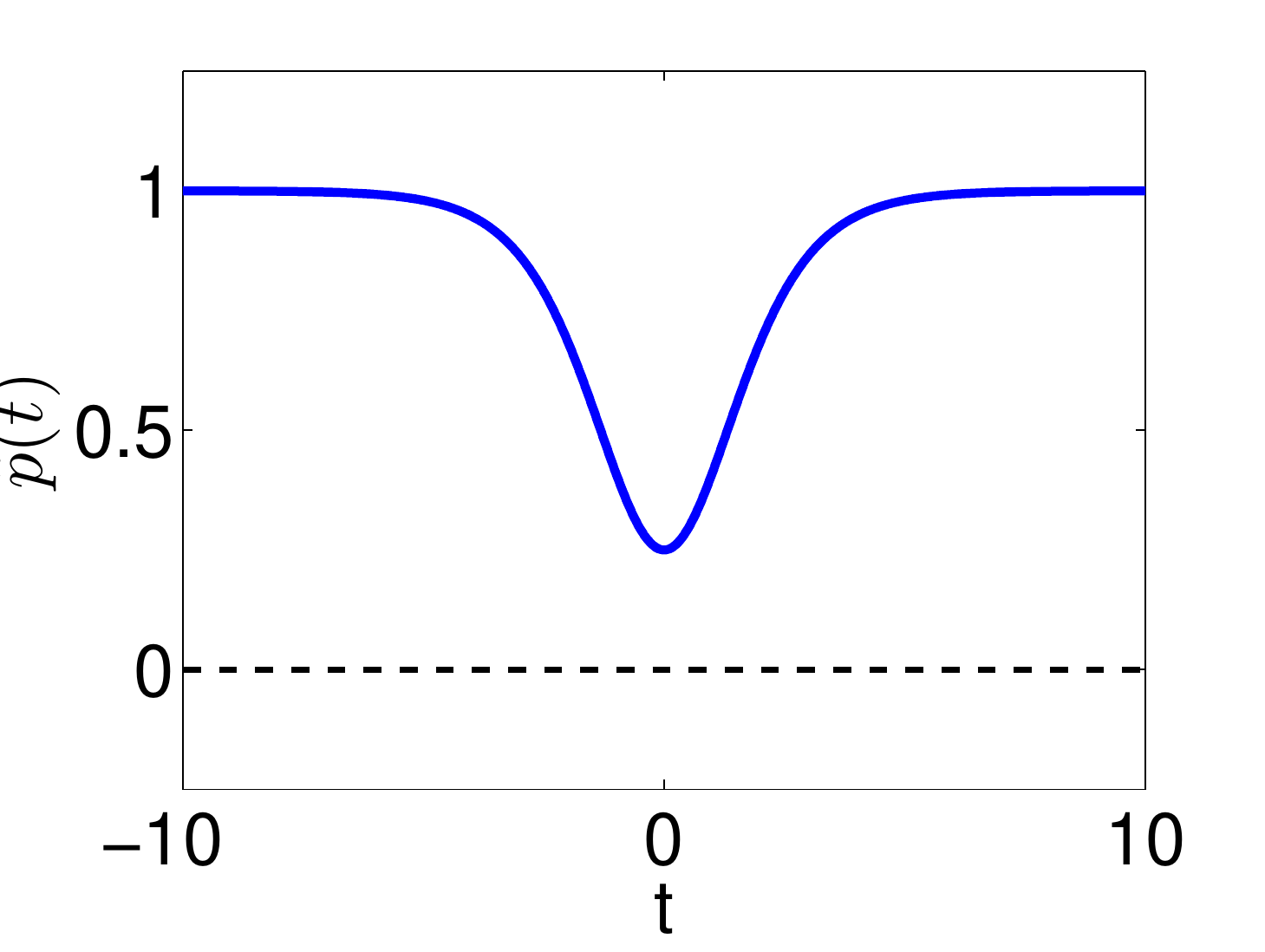}}
       \hfill
        \subcaptionbox{\label{SN nonlin drift TP}}[0.45\linewidth]
                {\includegraphics[scale = 0.3]{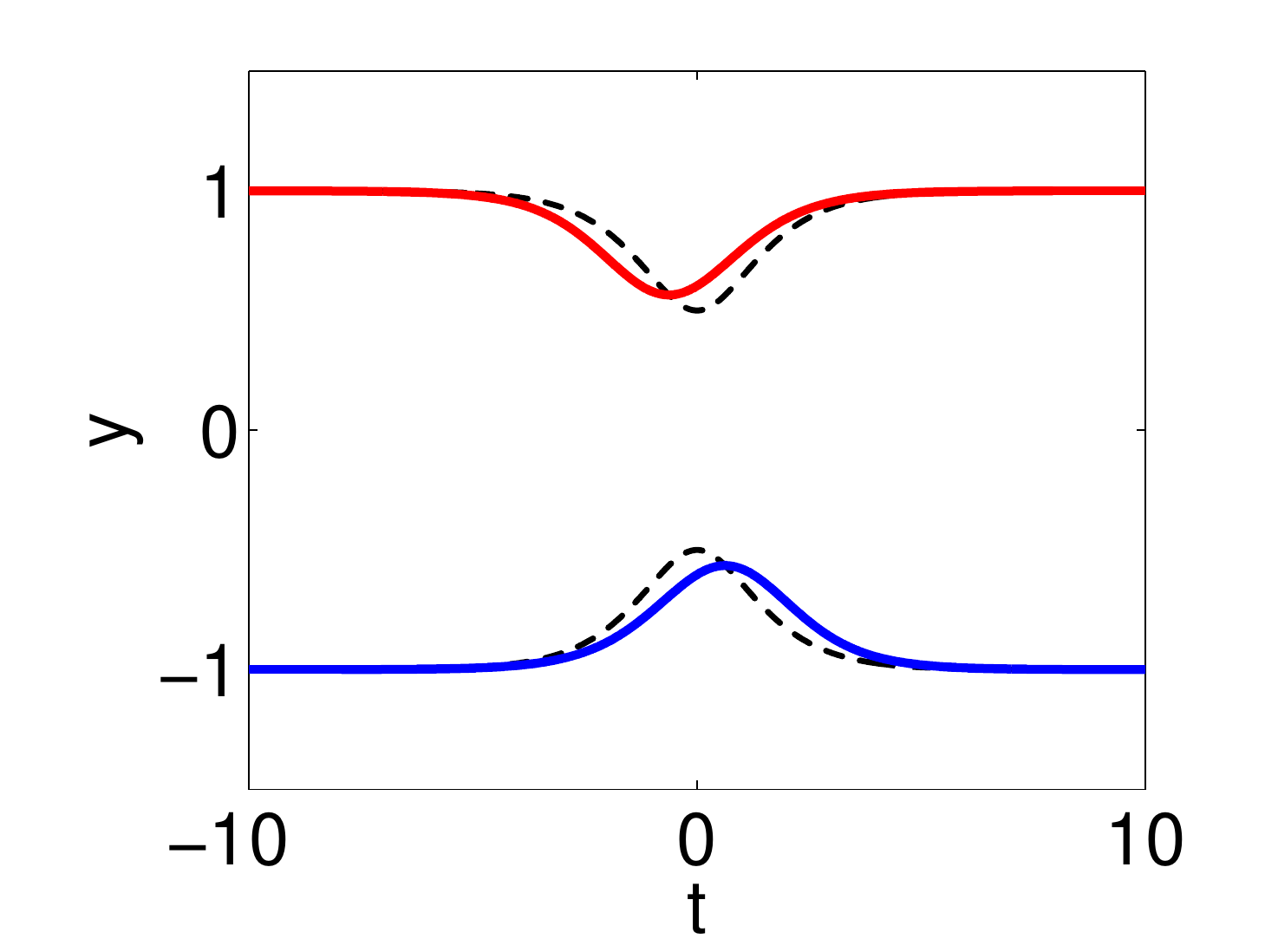}}
        \caption[Transformed nonlinear drift over same range as linear drift for saddle-node bifurcation.]{(a) Transformed nonlinear drift \eqref{transformed drift} for $\epsilon = 1/3$ (blue) such that range of values are the same as the linear example and so the saddle-node bifurcation (black dashed line) is not crossed. (b) Time profile for trajectories of a realisation that either starts at the stable quasi-static euilibrium (blue) or finishes at the unstable quasi-static equilibrium (red) for the system \eqref{transformed ODE}--\eqref{transformed drift}. Lower (upper) black dahed curves stable (unstable) quasi-static equilibrium. Parameters: $p_{0} = 1$, $\lambda_{\max} = 3$.}\label{SN nonlinear drift}
\end{figure}

However, to have some resemblance with the linear drift example, we will consider a drift that operates on the same range, namely $p(t) \in [0.25,1]$  for all $t$ and thus, does not cause the system to pass through the saddle-node bifurcation. This corresponds to a choice of $\epsilon = 1/3$ for which the parameters $p_{0} = 1$ and $\lambda_{\max} = 3$ remain the same, see Figure \ref{SN nonlin drift}. The time profile for this choice of drift is given in Figure \ref{SN nonlin drift TP}. The lower (upper) black dashed line represents the stable (unstable) quasi-static equilibrium. Furthermore, notice that the trajectory for a realisation starting at the stable quasi-static equilibrium (blue) does not track exactly this equilibrium.  

\begin{figure*}
\centering
\includegraphics[scale=0.35]{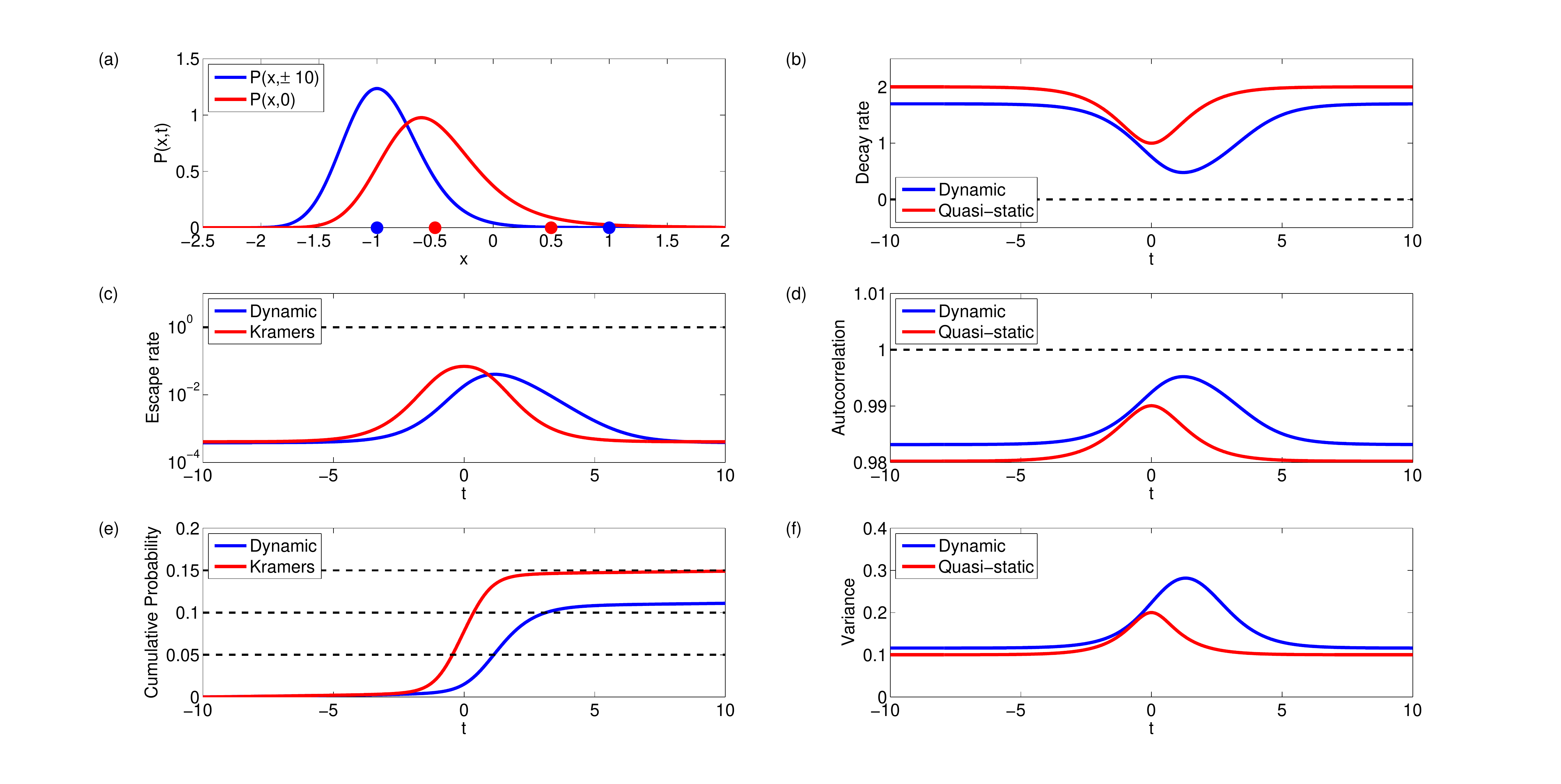}
\caption[Analysis of nonlinear drift towards and away from a saddle-node bifurcation.]{Results of a nonlinear drift towards and away from a saddle-node bifurcation according to the SDE \eqref{SN SDE} and the corresponding Fokker-Planck equation \eqref{Fokker Planck}. Nonlinear drift \eqref{transformed drift}, $p_{0} = 1$, $\epsilon = 1/3$, $D = 0.2$, $\Delta x = 0.05$, $\Delta t = 0.01$. (a) Stationary and final probability density $P_{s}(x) = P(x,\pm10)$ (blue) and density at $P(x,0)$ (red), \change{ fixed points indicated by dots on $x$-axis}. Decay rate (b), lag-1 autocorrelation (d) and variance (f) for both nonlinear dynamic (numerical calculation, blue) and linear quasi-static (linearised Ornstein-Uhlenbeck, red). Escape rate (c) and Cumulative probability of escape (e) for dynamic (blue) and Kramers' (red).}
\label{SN eg2b labels}
\end{figure*}

Using the SDE \eqref{sec:SNexamples SDE} with our nonlinear drift \eqref{transformed drift} we can calculate the time evolution of the probability density from the Fokker-Planck equation \eqref{Fokker Planck}. We can then determine the early-warning indicators and escape rates, which are all presented in Figure \ref{SN eg2b labels}.

Panel (a) displays the probability densities $P(x,t)$ at three different time moments. The blue density represents both the stationary density $P_{s}(x)$ and is identical to the normalised density at $t = T_{\mathrm{end}}$. The density in red represents when the system is closest to the saddle-node bifurcation. The fixed points at the various time steps are represented by the dots on the $x$-axis. Similar to the linear drift example, the red density has spread out due to the potential well flattening and has a larger tail developing towards the unstable quasi-equilibrium. 

Kramers' escape rate is plotted in red in panel (c), which assumes the system is stationary and therefore the escape rate is symmetrical about $t = 0$ due to the symmetry in the drift. Whereas, the dynamic escape rate is not symmetrical and instead appears to contain a lag. Initially, the dynamic and Kramers' give a good match as with the linear drift example. However, as the drift increases and the bifurcation is approached, the dynamic escape rate does not give an instant response. Therefore, the peak of the escape (slightly lower than Kramers' peak) occurs at about $t = 1$ before reducing back towards the base level. The cumulative probability of escape, panel (e), demonstrates once again that Kramers' ($\sim15\%$) overs estimates the numerically calculated cumulative probability ($\sim11\%$). The cumulative probability takes the shape of a step function, indicating that there is a given period of time for which there is significant escape. Moreover, the cumulative probability as well as the escape rate tells us that Kramers' approximates more escape but for a shorter period of time compared with the dynamic escape. 

The decay rate and early-warning indicators: lag-1 autocorrelation and variance are presented in panels (b), (d) and (f) respectively. The decay rate of the linearised potential well (red) decreases as the saddle-node is approached and then increases at the same rate as we move away, which is as expected. Correspondingly the linear quasi-static autocorrelation increases as the bifurcation is approached and is symmetric about $t = 0$. However, the nonlinear dynamic lag-1 autocorrelation does not increase until later than the linear quasi-static. In addition, the autocorrelation keeps increasing even when moving away from the saddle-node. The nonlinear dynamic variance displays similar behaviour, which starts increasing on the approach to the bifurcation but has a peak after $t = 0$. Though the initial increase aligns more with the linear quasi-static variance than for the case of autocorrelation. The nonlinear dynamic lag-1 autocorrelation and variance both reiterate the belief that there is a lag in the system. It is therefore debatable for an example of bifurcation-induced tipping with certain nonlinear drifts whether the autocorrelation and variance give sufficient warning of the tipping event. 

The differences between the values of the decay rate, autocorrelation and variance is again believed to be linked to the linearisation and the drift in the system. We found previously that a large noise level increases the nonlinear dynamic autocorrelation and variance compared with the linear quasi-static. This explains the initial and final values when the system is close to stationary. Furthermore, increasing the speed of the drift brings the nonlinear dynamic autocorrelation and variance down in relation to the linear quasi-static, which potentially explains the comparison as the saddle-node is approached. Though there are further complications due to the apparent lag in the system and quick turnaround of the drift.

In the following section we will consider one of the recognised policy relevant climate tipping points the Indian monsoon \citep{lenton2008tipping}. In particular we will consider a model that was developed by \citet{zickfeld2005indian} and use it to determine if the early-warning indicators can help detect the `switching off' of the monsoon.

\section{Case Study: Indian Summer Monsoon}
\label{sec: ISM}

The Indian summer monsoon season tends to start around May/June and ends around September/October time \citep{zickfeld2004modeling}. \citet{zickfeld2004modeling} developed a reduced ODE based model (with no spatial resolution) that is designed to capture the key mechanisms of the Indian monsoon. It was demonstrated by \citet{zickfeld2005indian} using this reduced model that either increasing the planetary albedo or reducing the current $CO_{2}$ levels (less likely to occur in the real-world) would cause the system to pass through a saddle-node bifurcation. This would result in the `switching off' of the monsoon and therefore leaving a far drier climate than the current one. We aim to replicate the behaviour of the model from \citet{zickfeld2004modeling} with further simplifications to reduce the model to a scalar ODE. \change{Adding noise and using the Fokker-Planck equation associated to this SDE will allow us to test the early-warning indicators. The early-warning indicators: increase of autocorrelation and variance would inform if any warning of this bifurcation-induced tipping event can be detected when the planetary albedo is increased.}

\subsection{Introduction of model}

\begin{figure}[h!]
\centering
\includegraphics[scale=0.3]{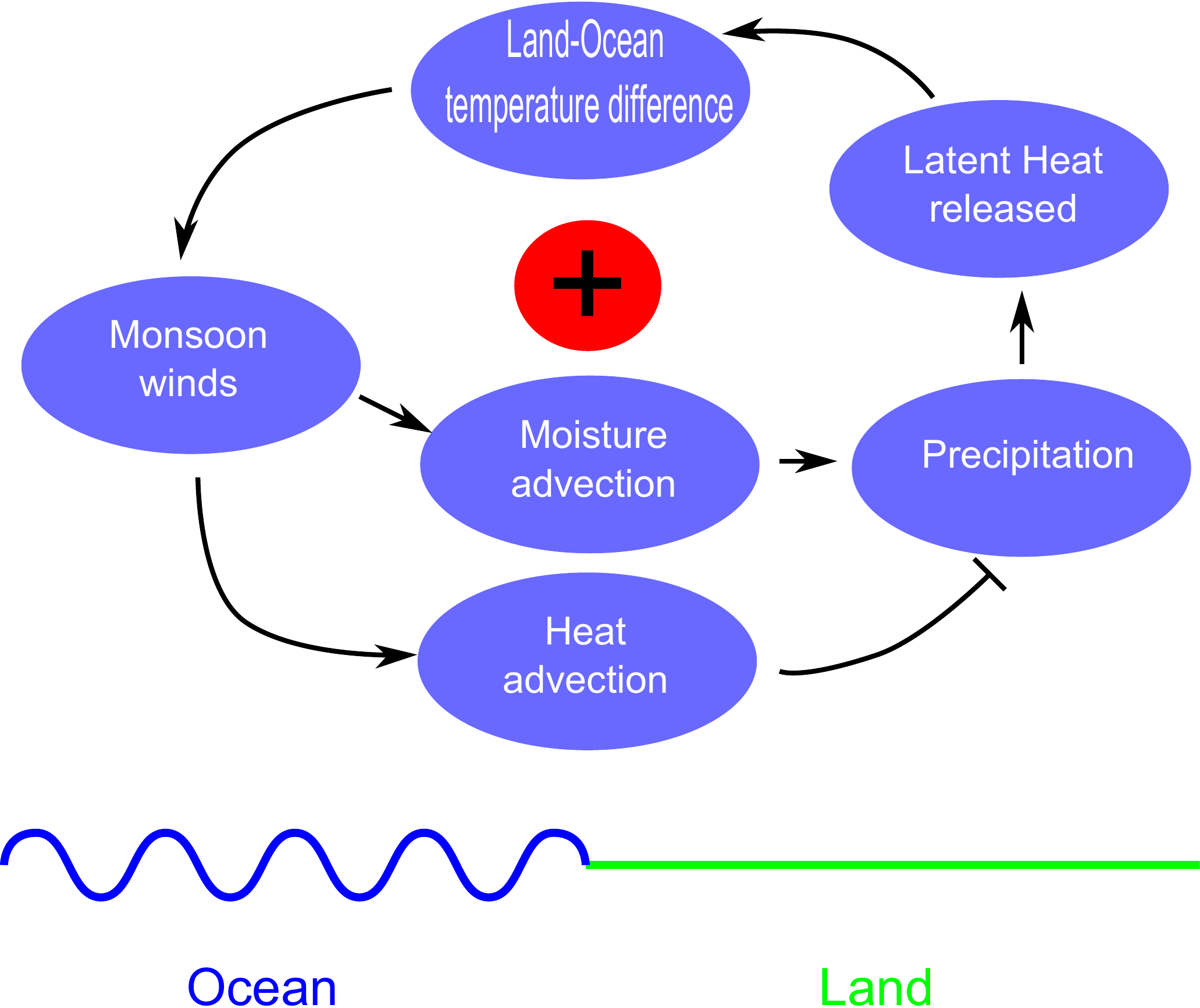}
\caption[Positive feedback loop for Indian summer monsoon.]{Positive feedback loop for Indian summer monsoon (based on \citet{levermann2009basic}).}
\label{Monsoon feedback}
\end{figure}

The driving force of the Indian summer monsoon is the positive feedback loop depicted in Figure \ref{Monsoon feedback}. During the winter months the prevailing winds over India are northeasterly (coming from the northeast) \citep{li2002winter}. This results in mainly dry winds coming from the Tibetan Plateau. However, as the temperature over land increases relative to the ocean, the winds are reversed such that they now come from the Indian Ocean \citep{zickfeld2004modeling}. The summer monsoon winds carry moisture from the Indian Ocean, which is then deposited as precipitation over India. This in turn releases latent heat and thus enhances the temperature difference between the land and ocean \citep{levermann2009basic}. The increase in the temperature gradient produces stronger winds coming off the ocean and so the positive feedback loop is formed.

We introduce the reduced model developed by \citet{zickfeld2004modeling} which captures the moisture-advection feedback as a key mechanism. The model treats India as a box, bounded by the Indian Ocean on three sides and the Tibetan Plateau to the north. A simplification that is made is to consider all four boundaries to be surrounded by the Indian Ocean. The atmosphere in the model is comprised of three layers and also contains two soil layers, where $w_{1}$ (top) and $w_{2}$ (bottom) represent the moisture in the respective layers. The soil layers represent two of the dynamic variables the other two are: $T_{a}$ the near surface air temperature and $q_{a}$ the specific humidity. \change{The specific humidity provides the moisture content in the air (ratio of water vapour mass to total mass of air parcel).} All other variables, for example precipitation, can be determined by these four dynamic variables. The ocean-land temperature difference $T_{a} - T_{oc}$ being a key quantity in many of the variables. 

The model is described by four ordinary differential equations:

\begin{align}
\dot{w_{1}} &= \dfrac{P - E - R}{f_{1}} + \dfrac{w_{2} - w_{1}}{\tau} \label{w1dot 2} \\
\dot{w_{2}} &= \dfrac{f_{1}(w_{1} - w_{2})}{f_{2}\tau} \label{w2dot 2} \\
\dot{q_{a}} &= \dfrac{E - P + A_{v}}{I_{q}} \label{qa dot 2} \\
\dot{T_{a}} &= \dfrac{\mathcal{L}(P - E) - F_{\uparrow}^{LW,TA} + F_{\downarrow}^{SL,TA}(1 - A_{\mathrm{sys}}) + A_{T}}{I_{T}} \label{Ta dot 2}
\end{align}

\noindent where the variables are summarised as follows:

\begin{itemize}
\item \textbf{Evaporation $(mm/s)$:} directly proportional to the soil moisture in the top layer, temperature difference between the ocean $T_{oc}$ and land $T_{a}$ and difference between saturated $q_{\mathrm{sat}}(T_{s})$ and specific humidity $q_{a}$.
\begin{align}
\nonumber
E&:= E(w_{1},q_{a},T_{a}) \\
&= Aw_{1}(T_{a}-T_{oc})(q_{\mathrm{sat}}(T_{s})-q_{a}) = w_{1}\tilde{E}(q_{a},T_{a})
\label{Evaporation eq}
\end{align}
\item \textbf{Precipitation $(mm/s)$:} directly proportional to the specific humidity.
\begin{equation}
P:= P(q_{a}) = Bq_{a}
\label{Precipitation eq}
\end{equation}
\item \textbf{Runoff $(mm/s)$:} directly proportional to soil moisture in top layer and the amount of precipitation fallen.
\begin{equation*}
R := R(w_{1},q_{a}) = Cw_{1}P(q_{a}) = w_{1}\tilde{R}(q_a)
\end{equation*}
\item \textbf{Moisture advection $(mm/s)$:} winds driven by the temperature gradient between land and ocean and are reversed above a monsoon turning height $h$. Below $h$ moisture is advected in proportional to $q_{oc}$ and above $h$ away proportional to $q_{a}$. 
\begin{equation}
A_{v} := A_{v}(q_{a},T_{a}) = G(T_{a}-T_{oc})(g_{1}q_{oc}-g_{2}q_{a})
\label{Moisture adv eq}
\end{equation}
\item \textbf{Net incoming solar radiation $(Kg/s^{3})$:} fraction of incoming solar radiation $I_{0}\cos\xi$ that is not reflected, determined by the planetary albedo $A_{\mathrm{sys}}$.
\begin{equation*}
F_{\downarrow}^{SL,TA}(1 - A_{\mathrm{sys}}) = I_{0}\cos\xi(1-A_{\mathrm{sys}})
\end{equation*}
\noindent where $\xi$ is the solar zenith angle - the measured angle of the sun's position from vertically above.
\item \textbf{Outgoing Long-wave radiation $(Kg/s^{3})$:} directly proportional to the temperature.
\begin{equation*}
F_{\uparrow}^{LW,TA}:= F_{\uparrow}^{LW,TA}(T_{a}) = HT_{a}+J
\end{equation*}
\item \textbf{Heat advection $(Kg/s^{3})$:} winds driven by temperature gradient, reversed above monsoon turning height $h$ determined by specific humidity, $q_{a}$. Cooler temperature is advected in for heights $z_{1}< h$ proportional to potential temperature above ocean $\theta_{oc}(T_{oc},z_{1})$ and advected away for heights  $z_{2} > h$ proportional to potential temperature above land $\theta_{a}(T_{a},z_{2})$. \change{The potential temperature at a height $z$ for a surface temperature $T_{s}$ and humidity $q_{s}$ is defined as $\theta_{s} = T_{s} - (\Gamma(T_{s},q_{s}) - \Gamma_{a})z$ where $\Gamma$ is the atmospheric lapse rate and $\Gamma_{a}$ the adiabatic lapse rate.}
\begin{align}
\nonumber
A_{T}&:= A_{T}(q_{a},T_{a}) \\
&= K(T_{a}-T_{oc})(\theta_{oc}(T_{oc},z_{1})-\theta_{a}(T_{a},z_{2})
\label{Heat adv eq}
\end{align}
\end{itemize}

\noindent The parameters $f_{1}$ and $f_{2}$ represent the field capacity of the upper and lower soil layers respectively. The soil moisture diffusion specific time is given by $\tau$ and $\mathcal{L}$ is the latent heat of evaporation. The remaining parameters $I_{q}$, $I_{T}$, $A-C$, $G-H$, $J$, $g_{1}$ and $g_{2}$ are all constants. 

\change{Equations \eqref{Evaporation eq}--\eqref{Heat adv eq} demonstrate how the positive feedback, depicted in Figure \ref{Monsoon feedback} features in the model \eqref{w1dot 2}--\eqref{Ta dot 2}. In the summer the temperature difference $T_{a} - T_{oc}$ increases, which causes an increase in evaporation \eqref{Evaporation eq} and moisture advection \eqref{Moisture adv eq} but has a negative impact on the heat advection \eqref{Heat adv eq}. The increase of evaporation and moisture advection increases the humidity \eqref{qa dot 2}. An increase of specific humidity leads to an increase of precipitation \eqref{Precipitation eq}. Consequently the atmospheric temperature increases \eqref{Ta dot 2} (despite small negative feedbacks) which boosts the land-ocean temperature difference and completes the positive feedback.}

\subsection{Steady state solutions}

We want to reduce the system \eqref{w1dot 2}--\eqref{Ta dot 2} from a four dimensional model down to a scalar ODE so that we use the Fokker-Planck equation to determine the early-warning indicators. We begin by removing the soil moisture feedback from the model. \change{The underlying mechanisms of the soil moisture in the two layers $w_{1}$, $w_{2}$ are fast such that they relax to the equilibrium rapidly compared to $q_{a}$ and $T_{a}$. We therefore set $w_{1}$, $w_{2}$ to equilibrium by setting \eqref{w1dot 2}--\eqref{w2dot 2} to zero.} Setting \eqref{w2dot 2} to equilibrium determines that the soil moisture is the same in both the upper and lower layers. Thus from \eqref{w1dot 2} we have:

\begin{equation*}
P(q_{a}) - w_{1}(\tilde{E}(q_{a},T_{a}) + \tilde{R}(q_{a})) = 0
\end{equation*}  

\noindent and hence the soil moisture in the two layers can be determined from

\begin{equation*}
w_{1} = w_{2} = \dfrac{P(q_{a})}{\tilde{E}(q_{a},T_{a}) + \tilde{R}(q_{a})}
\end{equation*}

\noindent using present day values for the specific humidity $q_{a}$ and temperature $T_{a}$. Removing the soil moisture feedback from the model reduces the system to two dimensions \eqref{qa dot 2}--\eqref{Ta dot 2}.

\begin{figure}[h!]
\centering
\includegraphics[scale=0.3]{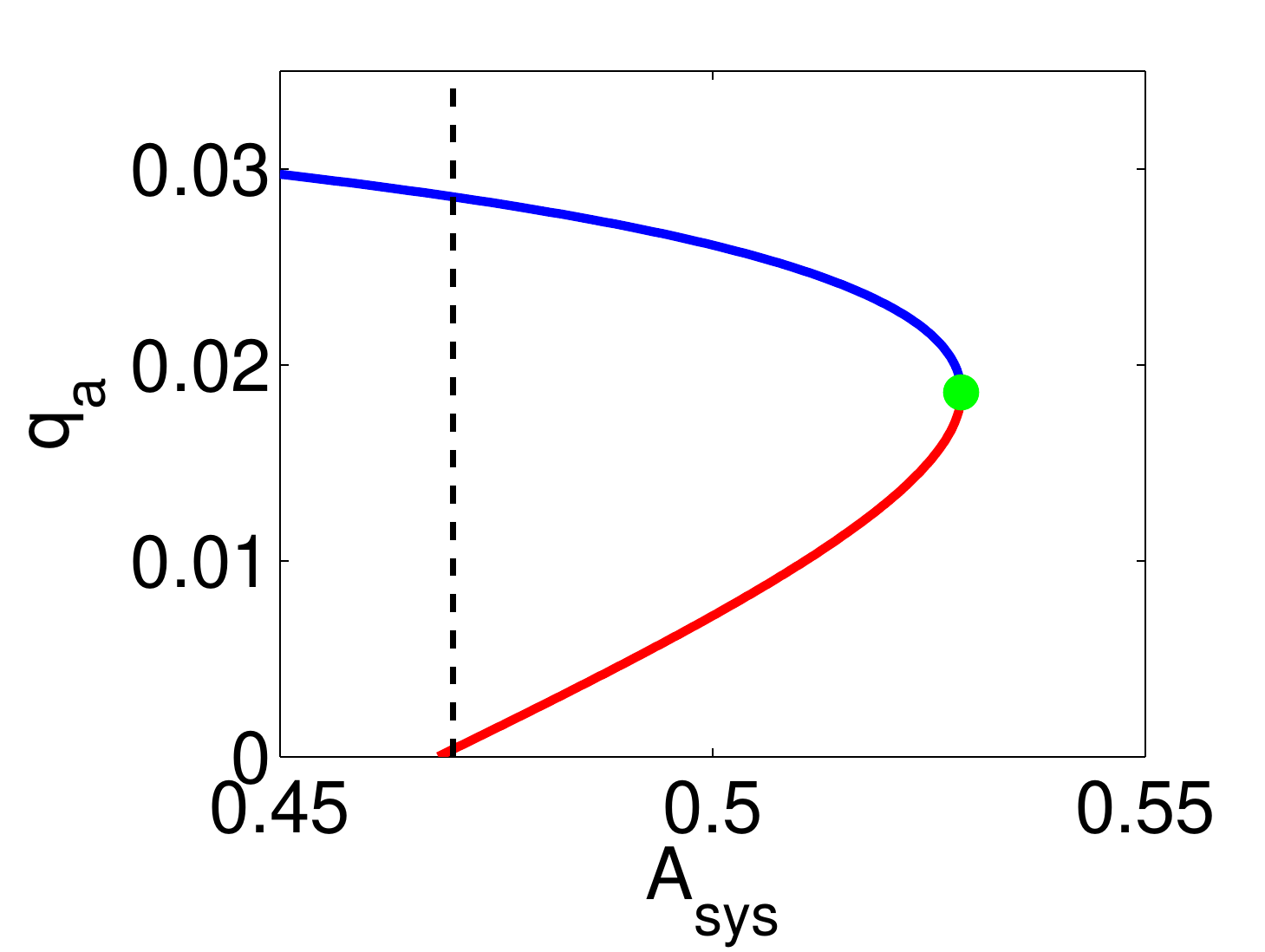}
\caption[Bifurcation diagram for Indian summer monsoon model.]{Bifurcation diagram for Indian summer monsoon model. Saddle-node bifurcation point indicated by green dot, upper branch is stable (blue) and lower branch unstable (red). Present day value of planetary albedo $A_{\mathrm{sys}}$ indicated by vertical black dashed line.}
\label{monsoon bif}
\end{figure}

\citet{zickfeld2005indian} identified two parameters, $CO_{2}$ concentration and the planetary albedo, which play a crucial role in the destabilising mechanism of the system and are influenced by human activities or subject to natural variation. Reducing $CO_{2}$ emissions, causes the system to pass through a saddle-node bifurcation (not shown). However, we will consider the planetary albedo $A_{\mathrm{sys}}$ at the top of the atmosphere. The planetary albedo represents the ratio of reflected to incoming solar radiation and can be affected by atmospheric aerosols and land cover conversion \citep{zickfeld2005indian}. We will therefore use the planetary albedo as our bifurcation parameter and focus on the specific humidity $q_{a}$ as the albedo is increased. Changes to the specific humidity will indicate whether the summer monsoon is in an `on' or `off' state by the amount of moisture in the air. We can ascertain the bifurcation diagram by setting \eqref{qa dot 2} to zero and derive a formula for the temperature $T_{a}$ dependent on the humidity $q_{a}$. Finally setting \eqref{Ta dot 2} to equilibrium allows us to determine the planetary albedo $A_{\mathrm{sys}}$ as a function of $q_{a}$. The bifurcation diagram is presented in Figure \ref{monsoon bif}.   

The present day value of the planetary albedo, $A_{\mathrm{sys}}$, is $0.47$ indicated by the vertical black dashed line. The blue branch is the family of stable equilibria and red is the unstable equilibrium branch, thus, the current specific humidity, $q_{a}$, is about $0.03$. A specific humidity of $0.03$ indicates moisture in the air and represents the monsoon being in an `on' state. In addition, if the planetary albedo was to change a small amount, for example, decreasing or increasing by $0.02$ the specific humidity would remain about the same. However, a further increase of the planetary albedo brings the system ever closer to the saddle-node bifurcation, indicated by the green dot. The system is then susceptible to tipping for sufficient noise in the system, which will cause the specific humidity to be close to zero indicating little moisture in the air and thus, the monsoon is in an `off' state.   

We will simulate this behaviour and observe if the early-warning indicators are present as the bifurcation is approached. Close to the saddle-node bifurcation the dynamics will be dominated by the slow variable, which is assumed to be the specific humidity $q_{a}$. We therefore reduce the two dimensional system \eqref{qa dot 2} -- \eqref{Ta dot 2} down to a single ODE by setting the atmospheric temperature $T_{a}$ to equilibrium. Setting \eqref{Ta dot 2} to zero and rearranging gives a quadratic expression for $T_{a}$ (we choose the physically relevant root) that depends on the specific humidity and planetary albedo. We add white noise to the dynamics of $q_{a}$, equation \eqref{qa dot 2}, to create the SDE:

\begin{align}
\nonumber
\mathrm{d}Q_{t} = \bigg[&\dfrac{E(Q_{t},A_{\mathrm{sys}}(t)) - P(Q_{t}) + A_{v}(Q_{t},A_{\mathrm{sys}}(t))}{I_{q}}\bigg]\mathrm{dt} \\
&+ \sqrt{2D}\mathrm{d}W_{t}
\label{qa SDE}
\end{align}   

\noindent where the planetary albedo is the driving force of the system, given by:

\begin{equation*}
A_{\mathrm{sys}}(t) = A_{\mathrm{sys},0} + \epsilon t
\end{equation*}

\noindent and $A_{\mathrm{sys},0} = 0.47$ is the present day value of the planetary albedo. The drift speed $\epsilon = 0.006$ is chosen such that over the span of ten decades the planetary albedo increases from it's present day value to just past the saddle-node bifurcation.

\subsection{Presence of early-warning indicators for simulated monsoon tipping}

\begin{figure*}
\centering
\includegraphics[scale=0.35]{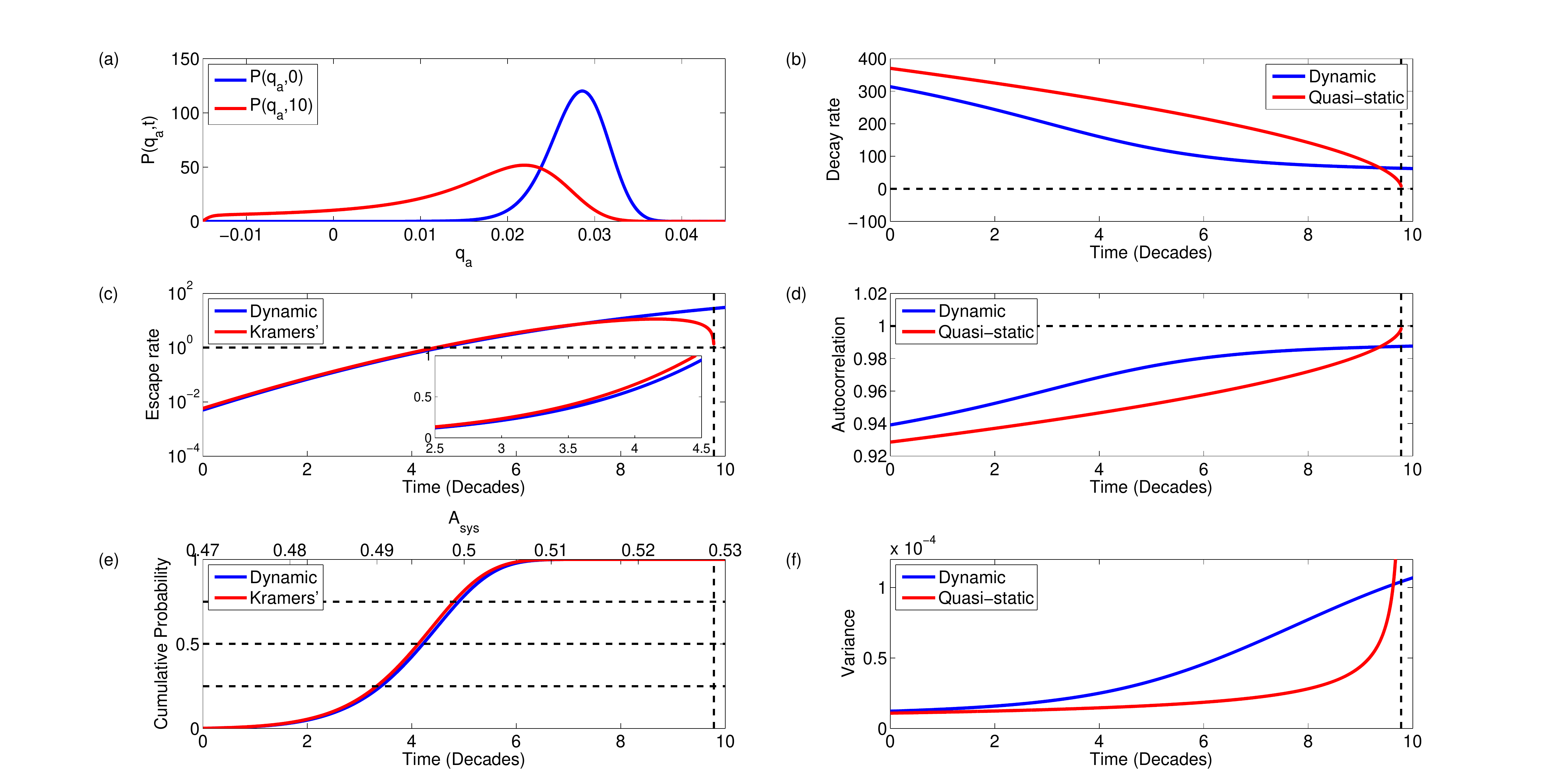}
\caption[Analysis of slow passage through a saddle-node bifurcation for an Indian summer monsoon model.]{Results of slow passage through a saddle-node bifurcation for an Indian summer monsoon model according to the SDE \eqref{qa SDE} and the corresponding Fokker-Planck equation \eqref{Fokker Planck}. Black vertical dashed lines indicate time saddle-node bifurcation is reached. Linear drift: $A_{\mathrm{sys}}(t) = A_{\mathrm{sys},0} + \epsilon t$, $A_{\mathrm{sys},0} = 0.47$, $\epsilon = 0.006$, $D = 0.004$, $\Delta x = 0.001$, $\Delta t = 0.0002$. (a) Stationary and final probability density $P_{s}(q_{a}) = P(q_{a},0)$ (blue) and density at $P(q_{a},10)$ (red). Decay rate (b), lag-1 autocorrelation (d) and variance (f) for both nonlinear dynamic (numerical calculation, blue) and linear quasi-static (linearised Ornstein-Uhlenbeck, red). Escape rate (c) and Cumulative probability of escape (e) for dynamic (blue) and Kramers' (red). Present day parameter values sourced from \citet{zickfeld2004modeling}, future projections $(\epsilon,D)$ chosen for illustration purposes.}
\label{monsoon new3}
\end{figure*}

We generate the probability density function $P(q_{a},t)$ using the Fokker-Planck equation \eqref{Fokker Planck} with the drift $f(q_{a},t)$ given by the expression inside the square brackets in \eqref{qa SDE}. We subsequently calculate the escape rate and the early-warning indicators: variance and autocorrelation using \eqref{sec:SN variance} and \eqref{autocorrelation SN sec} respectively. The results are presented in Figure \ref{monsoon new3}.

Panel (a) depicts the normalised stationary density $P_{s}(q_{a}) = P(q_{a},0)$ (blue) and the end density $P(q_{a},T_{\mathrm{end}})$ (red). The stationary density is approximately of the form of a normal distribution centred about the stable equilibrium at about $0.03$. As the density evolves over time it begins to widen as the potential well flattens. For $t = T_{\mathrm{end}}$ no equilibria exist though the density still has a small peak centred at about $q_{a} = 0.02$ but a very wide tail with some escape. A peak still exists as the bifurcation has only just been crossed and so the potential is still relatively flat. However, if we were to go sufficiently past the bifurcation the density would build up on the left boundary.   

The escape rate and cumulative probability of escape are represented in panels (c) and (e) respectively. There is a good agreement between Kramers' (red) and the dynamic (blue) escape rates for $t \leq 8$ decades. Over this period Kramers' escape rate provides a slight overestimation as can be seen by the insert in panel (c). However, Kramers' escape rate fails to give an accurate representation of the escape rate close to or past the saddle-node bifurcation. This is because \change{one assumption behind the validity of Kramers' escape formula is that the potential landscape is stationary, whereas, we use a slow drift which explains the initial discrepancy far away from the saddle-node bifurcation. Additionally Kramers' escape formula assumes that the noise level $D$ is much smaller than the well depth. Though close to a saddle-node bifurcation the noise level is always relatively large because the well is shallow.} The cumulative probability of escape provides an estimation for when the tipping occurs. \change{After approximately $40$ years or at $A_{\mathrm{sys}} = 0.495$ a realisation has escaped with $50\%$ probability, which} is a relatively long time before the bifurcation is reached at about $98$ years. 

We examine the decay rate (b) and the early-warning indicators; lag-1 autocorrelation (d) and variance (f), to see if they offer any forewarning of the approach of a bifurcation-induced tipping event. The nonlinear dynamic decay rate (blue) decreases as the bifurcation is approached and thus, correspondingly the nonlinear dynamic autocorrelation (blue) increases. Combining this with the nonlinear dynamic variance (blue), which also increases we can state that the early-warning indicators do offer forewarning of a bifurcation-induced tipping event. The linear quasi-static approximations (red) of the early-warning indicators produce similar behaviour, namely the autocorrelation and variance increases. However, the linear quasi-static approximations underestimate the nonlinear dynamic indicators, which, as discussed in Section \ref{sec: SN examples}, can be attributed to the relatively large noise level. Furthermore, the shape between the nonlinear dynamic and linear quasi-static variances is different, this can largely be attributed to the domain size $[-0.015,0.045]$ and therefore restricts the nonlinear dynamic variance. 

\begin{figure}[h!]
        \centering
        \subcaptionbox{\label{Monsoon prob D}}[0.45\linewidth]
                {\includegraphics[scale = 0.3]{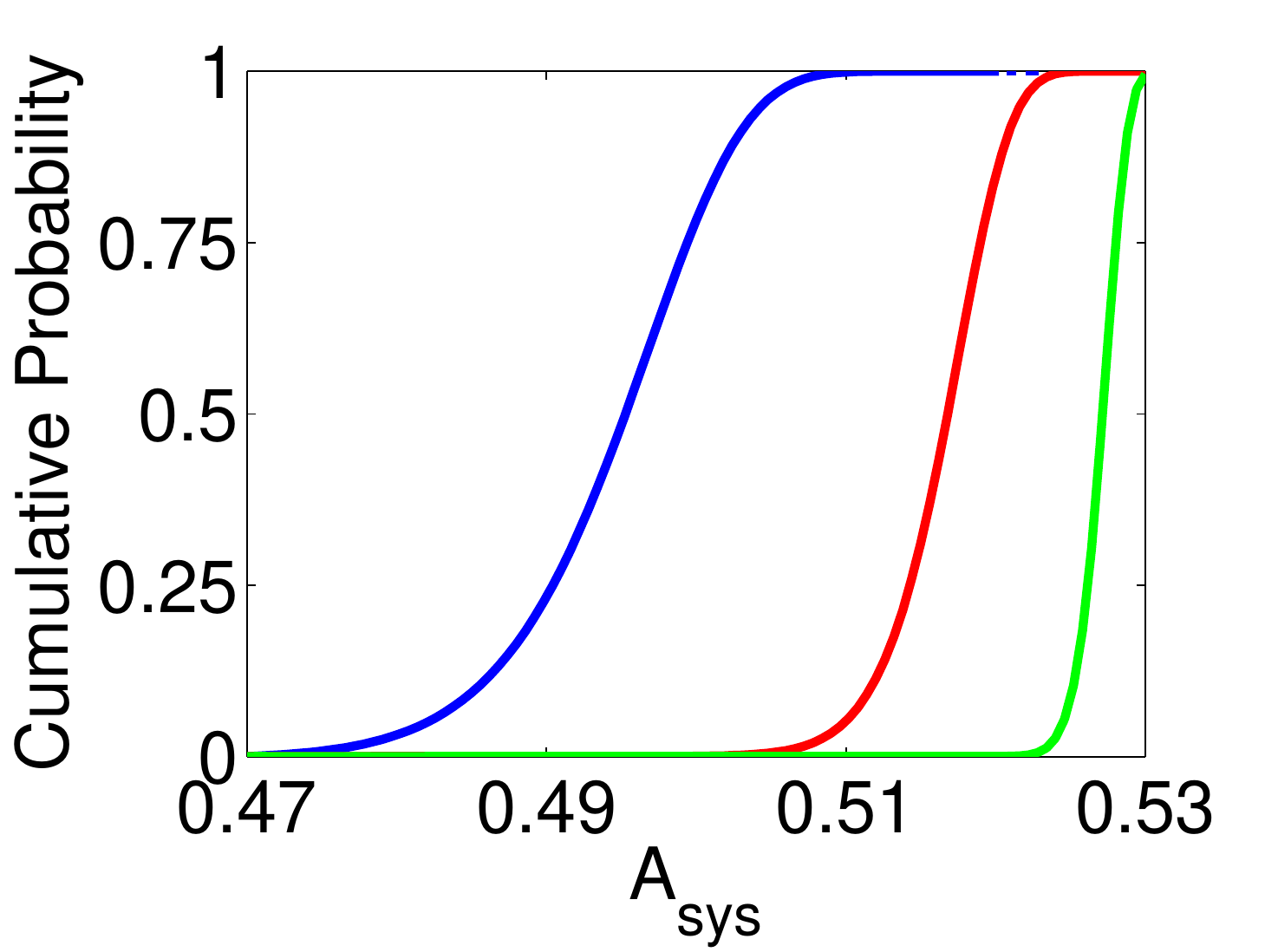}}
       \hfill
        \subcaptionbox{\label{Monsoon prob epsilon}}[0.45\linewidth]
                {\includegraphics[scale = 0.3]{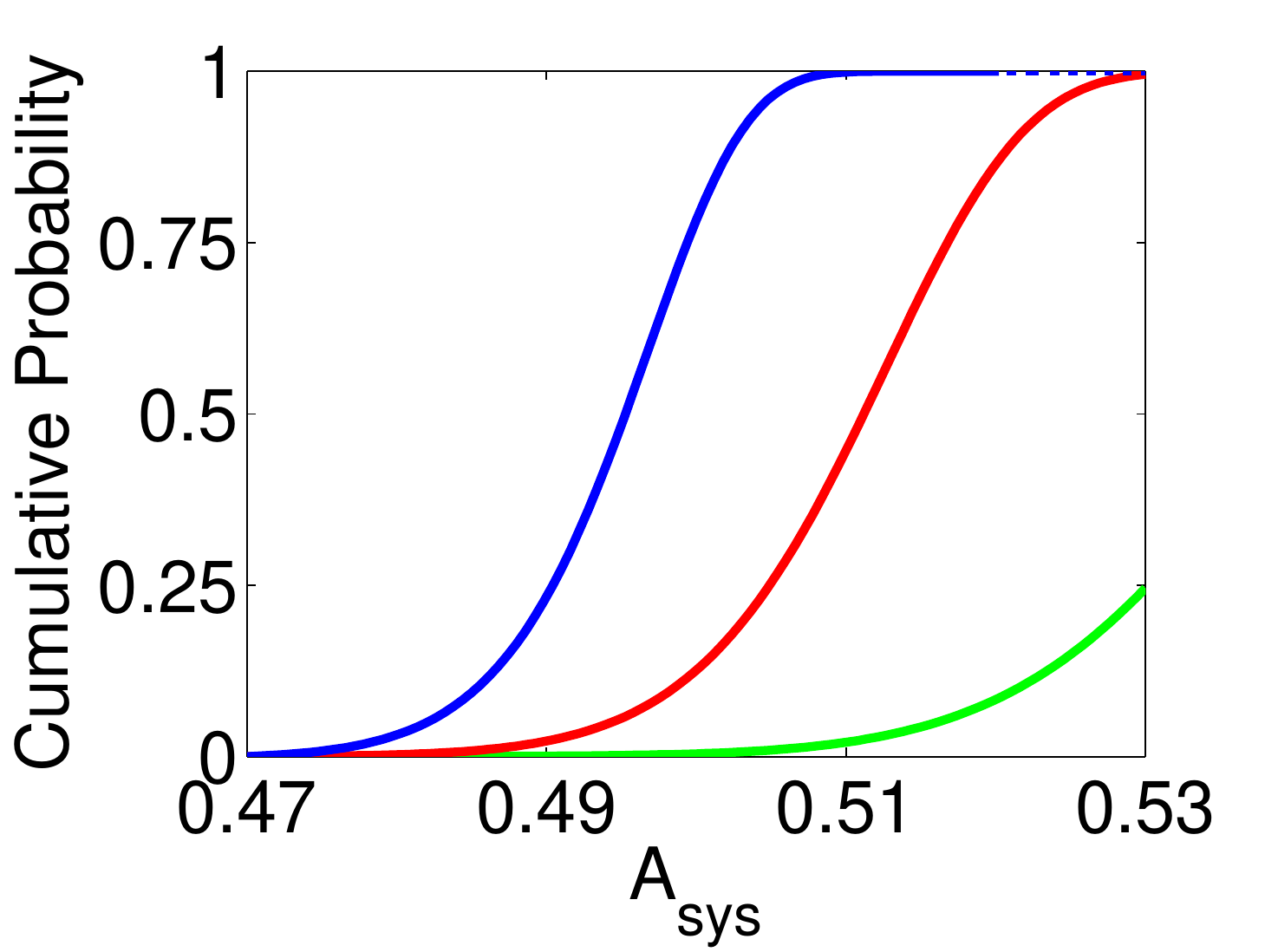}}
        ~ 
        \caption[Cumulative probability of escape of monsoon model.]{Comparing cumulative probability of escape for monsoon model with different noise levels (a) and drift speeds corresponding to how quickly planetary albedo $A_{\mathrm{sys}}$ is increased (b). (a) Noise levels used $D = 0.004$ (blue), $D = 0.0012$ (red), $D = 0.0002$ (green), for $A_{\mathrm{sys}}$ increased over $100$ year period. (b) Planetary albedo $A_{\mathrm{sys}}$ changed over $100$ years (blue), $10$ years (red) and $1$ year period (green) with noise level $D = 0.004$.}\label{Monsoon cum. prob}
\end{figure}

In Figure \ref{Monsoon cum. prob} we analyse how the noise level and drift speed (how quickly planetary albedo is increased) affect the cumulative probability and thus timing of the tipping. For a slowing changing planetary albedo and relatively large noise we found that the tipping occurred at about $A_{\mathrm{sys}} = 0.495$ (blue curves). This equates to a long time before the bifurcation point at about $A_{\mathrm{sys}} = 0.529$, we first investigate how the noise level affects the cumulative probability of escape.   

In Figure \ref{Monsoon prob D} we keep the drift speed fixed and change the noise level, where the largest noise level is given in blue and the smallest noise level in green. We see that as expected a smaller noise level causes a sharper and later transition in the cumulative probability. This is because for a small noise level (green) a realisation has only a `small window' close to the bifurcation to escape. Whereas, with a large noise level realisations are more likely to escape further from the bifurcation. Therefore these realisations have a `wider window' in when it is possible to escape leading to a shallower increase in the cumulative probability.

Figure \ref{Monsoon prob epsilon} addresses how the speed of albedo change affects the cumulative probability for a fixed noise level. The blue curve equates to changing the albedo over the given range in the space of 100 years. Increasing the planetary albedo in a 10 year period is given in red and over just a year is given in green. At first glance Figure \ref{Monsoon prob epsilon} can be a little misleading as the cumulative probability is plotted over the planetary albedo as opposed to time. Therefore it should be noted that the quicker the transition of the planetary albedo the sharper the transition in the cumulative probability in time. For a slow transition with relatively large noise (blue) there is a $50\%$ chance the system has tipped by about $A_{\mathrm{sys}} = 0.495$. Increasing the planetary albedo over a 10 year period (red) $50\%$ escape is reached closer to the bifurcation at about $A_{\mathrm{sys}} = 0.51$. Though for an extremely fast change in the albedo over only a 1 year period (green) we find that the cumulative probability of escape at $A_{\mathrm{sys}} = 0.53$ is only $0.25$. Therefore even having passed the bifurcation most realisations have not escaped yet. 

This behaviour is similar \change{to that for a rapid shift $(\epsilon)$ the escape is delayed (same as small noise, see saddle-node normal form Section \ref{subsec: SN nlin drift}}). In Section \ref{subsec: SN nlin drift} we observed that it was possible to have nonlinear drifts that would take the system briefly past the bifurcation and then back without causing a transition. In this example we have just passed the bifurcation and yet to make the transition, though continued increase of the albedo will cause a tipping. A question often raised by scientists but as of writing not addressed is, if a system is known to be just beyond the tipping point but has not yet tipped what needs to be done to stop the tipping all together? This would normally require reversing the direction of the bifurcation parameter at a given speed. The nonlinear shift (Section \ref{subsec: SN nlin drift}), motivated from rate-induced tipping is an example of this occurring but does not inform the minimum drift required.    

\bibliography{SN_bib}

\end{document}